\documentclass[
final,
a4paper,
11pt,
]{article}
\usepackage{preamble}
\usepackage{tikzsetting}

\makeatletter
\let\ftype@table\ftype@figure
\makeatother

\title{Artificial compressibility methods for the incompressible Navier--Stokes
equations using lowest-order face-based schemes on polytopal meshes}
\date{}
\author{Riccardo Milani\footnotemark[1]\textsuperscript{ , }\footnotemark[2]
\and J\'er\^ome Bonelle\thanks{EDF R\&D, 6 Quai Watier, 78400 Chatou, France}
\and Alexandre Ern\thanks{CERMICS, Ecole des Ponts, 77455 Marne-la-Vall\'{e}e 2, France, 
and  Inria, 2 rue Simone Iff, 75589 Paris, France}
}

\begin{document}
\maketitle

\begin{abstract}
We investigate artificial compressibility (AC) techniques for the time
discretization of the incompressible Navier--Stokes equations. The space
discretization is based on a lowest-order face-based scheme supporting
polytopal meshes, namely discrete velocities are attached to the mesh faces and
cells, whereas discrete pressures are attached to the mesh cells. This
face-based scheme can be embedded into the framework of hybrid mixed mimetic
schemes and gradient schemes, and has close links to the lowest-order version
of hybrid high-order methods devised for the steady incompressible
Navier--Stokes equations. The AC time-stepping uncouples at each time step the
velocity update from the pressure update. The performances of this approach are
compared against those of the more traditional monolithic approach which
maintains the velocity-pressure coupling at each time step. We consider both
first-order and second-order time schemes and either an implicit or an
explicit treatment of the nonlinear convection term. We investigate
numerically the CFL stability restriction resulting from an explicit treatment,
both on Cartesian and polytopal meshes. Finally, numerical tests on large 3D
polytopal meshes highlight the efficiency of the AC approach and the benefits
of using second-order schemes whenever accurate discrete solutions are to be
attained. 


\end{abstract}

\section{Introduction}

The goal of this work is to study the accuracy and efficiency of artificial
compressibility techniques for the time discretization of the incompressible
Navier--Stokes equations. These equations are encountered in a wide range of
industrial applications, ranging from aeronautics to the simulations of flows in
micro-fractures, to cite two salient examples. In this work, the space discretization
is realized by means of lowest-order face-based schemes supporting classical
(simplicial or Cartesian) as well as general (polytopal) meshes. The support of
polytopal meshes is motivated by the applications,  especially in the industrial
context. Indeed, such meshes often alleviate substantially the burden of mesh
generation resulting from the complexity of the geometry (as, e.g., the shape of
industrial equipments) and/or the heterogeneity of the physical properties (requiring
local mesh refinements producing hanging nodes).

The unsteady incompressible Navier--Stokes equations read as follows:
\begin{equation}
  \label{eq:NS}
  \left\{
    \begin{aligned}
      \ddt{\Vex}-\Visc\lapl\Vex +(\Vex\cdot\grad)\Vex +\grad\Pex ={}&\Foex\,,\\
      \dive\Vex                                                  ={}&0\,,
    \end{aligned}
  \right. \qquad \text{in } \Dom\times(0,T)\,,
\end{equation}
together with suitable boundary and initial conditions on  the velocity. Here,
$\Dom\subset\Real^{\dmn}$, $\dmn=2,3$, is an open, polytopal, bounded, Lipschitz
domain and $T>0$ denotes the observation time. The unknowns are the velocity $\Vex$
and the pressure $\Pex$ (by convention, vector-valued quantities are underlined). The
problem data are the viscosity $\Visc>0$, the mass density is normalized to unity,
and the body force is denoted by $\Foex$. For simplicity, we consider smooth
solutions to \eqref{eq:NS}, and we enforce Dirichlet boundary conditions over the
whole boundary $\pOm$ at all times. Recall that the first equation in~\eqref{eq:NS}
expresses the momentum balance in the flow, and the second equation the mass
conservation. 

Let us briefly describe the time-marching schemes considered in this work without
introducing any space discretization yet. These schemes are either of monolithic
form, thereby requiring to solve a saddle point problem at each time step, or they
uncouple the velocity and the pressure at each time step by means of an artificial
compressibility technique. Let $\deltat>0$ be the time step (taken to be constant for
simplicity) and let $(\Vex^n,\Pex^n)$ denote the approximate solution at  the
discrete time node $t^n\coloneqq n\deltat$ for all $n\geq1$. The monolithic approach
is the traditional time-marching scheme for the incompressible  Navier--Stokes
equations and it reads as follows:  For all $n\geq1$, given $\Vex\nm$ from the
initial condition or the previous time node, find $(\Vex^n,\Pex^n)$ by solving
the saddle point problem in $\Dom$,
\begin{equation}
  \label{eq:ust_mono}
  \left\{
    \begin{aligned}
      \dddtm{\Vex}-\Visc\lapl\Vex^n +(\Vex^n\cdot\grad)\Vex^n +\grad\Pex^n ={}&\Foex^n\coloneqq\Foex(t^n)\,,\\
      \dive\Vex^n ={}&0\,.
    \end{aligned}
  \right.
\end{equation}
The monolithic approach is well-known and is  considered here as the reference
time-stepping scheme. The only approximation introduced in \eqref{eq:ust_mono} with
respect to \eqref{eq:NS} is the discretization of the velocity time-derivative (here,
by means of the implicit Euler  scheme to fix the ideas). In contrast, the artificial
compressibility (AC) approach introduces an additional approximation at each time
step in that the velocity is first updated and the pressure is then corrected. In its
first-order form, the scheme reads as follows: For all $n\geq1$,  given
$(\Vex\nm,\Pex\nm)$ from the initial condition or the previous time node, find
$(\Vex^n,\Pex^n)$ by solving a parabolic-like problem on the velocity and then
updating the pressure:
\begin{subequations}
  \label{eq:ust_ac}
  \begin{align}
    \dddtm{\Vex}-\Visc(\lapl\Vex^n + \ACprm\graddiv\Vex^n) +(\Vex^n\cdot\grad)\Vex^n ={}&\Foex^n-\grad\Pex\nm\,,\label{eq:ust_ac-a}\\
    \Pex^n={}&\Pex\nm - \Visc\ACprm\dive\Vex^n\,,\label{eq:ust_ac-b}
  \end{align}
\end{subequations}
where $\ACprm>0$ is a user-defined adimensional parameter. The AC method was
introduced in the late sixties \cite{Chor68,Tema69a} in the Western literature,
although the seminal ideas can be traced independently in the Russian literature as
well \cite{VlKY66,Lady69,Yane71}. For a recent analysis of the method, the reader is
referred to \cite{GuMi15}. Notice the appearance of the additional grad-div term
in~\eqref{eq:ust_ac-a}, whereas \eqref{eq:ust_ac-b} shows that the time discrete 
velocity field $\Vex^n$ is no longer divergence-free. Notice also that an
approximation to the initial pressure is needed.  Furthermore, we observe that
higher-order versions of the AC scheme are available; see \cite{GuMi15,GuMi19}. 

The computational study
performed in this work considers both the monolithic and the AC approaches,
either in the above form which is first-order accurate in time, or in second-order
form (as described below). Moreover, we examine either the implicit treatment of the
convection term (as above) leading to a nonlinear problem to be solved at each time
step, or a (semi-)explicit treatment leading to a linear problem  to be solved at
each time step. The natural choice when using the AC scheme is to resort to
an explicit treatment since the goal is to reduce as much as possible the 
computational cost per time step. However, for completeness, we also
consider the implicit treatment of the convection term for 
the first-order AC schemes, but not for the second-order ones. 
Indeed, in this latter case, there are two nonlinear problems to be solved at each
time step, which is deemed to be too expensive.

Concerning the space discretization, we focus on schemes supporting polytopal meshes.
Numerous possibilities are available from the literature. Concerning high-order
schemes, we mention (without being exhaustive) discontinuous Galerkin (dG) methods
\cite{BaRe97,CKSS02,HaLa02}, hybridizable dG (HDG) methods
\cite{CoNP10,NPCo11,JePS14,LeSc16}, hybrid high-order (HHO) methods
\cite{ABDP15,DELS16,DPKr18,BDPD19}, weak Galerkin (WG) methods
\cite{MuWY15}, virtual element methods (VEM) \cite{BVLV17,BVLV18,GaMS18},
and nonconforming VEM \cite{CaGM16}. As efficient as the above high-order
methods can be, it is often preferred in an industrial context to use lowest-order
schemes owing to their relative ease of implementation and the prominence of legacy
codes based on these techniques. Moreover, since the maximal order of accuracy in
time considered here is two, it is reasonable to focus on lowest-order space
discretization schemes. Examples of such schemes for the (Navier--)Stokes equations
include  discrete duality finite volume (DDFV) schemes
\cite{Delc07,KrMa12,BoKN15,GoKL19}, mixed finite volume (MFV) schemes
\cite{DrEy09}, mimetic finite difference (MFD) schemes \cite{BGLM09},
compatible discrete operator (CDO) schemes \cite{BoEr15,Bone14}, and gradient
schemes \cite{Fero15,DrEF15,EyFG18}. Unifying frameworks bridging a large class
of lowest-order schemes exist, such as the hybrid mixed mimetic (HMM) framework from
\cite{DEGH13} and the setting of gradient schemes \cite{DEGH13}. In this
work, we focus more specifically on CDO schemes. This is a mild restriction, and we
expect that most of our conclusions can be carried out to other lowest-order schemes
owing to the above-mentioned unifying frameworks. Among the various CDO schemes, we
focus on face-based CDO (CDO-Fb) schemes in which the degrees of freedom (DoFs) for
the velocity are vector-valued and attached to the mesh faces and cells, whereas the
pressure DoFs are scalar-valued and attached to the mesh cells only. We refer the
reader to \cite{BoEr14} for the analysis of vertex-based and cell-based CDO
schemes for elliptic PDEs and to \cite{Cant16,CaEr17} for edge-based CDO
schemes applied to circulations (i.e., differential forms of order one). The CDO
schemes devised in \cite{BoEr15} are different since they introduce
the vorticity as an additional unknown. As discussed in
\cite[Sect.~8.3]{Bone14},  CDO-Fb schemes are derived for cell-based schemes by
a hybridization procedure of the flux unknown considered in cell-based schemes.
CDO-Fb schemes can be bridged to lowest-order HHO, HFV, HMM, and gradient schemes,
and have  been evaluated numerically to approximate the steady incompressible
Navier--Stokes equations in \cite{BoEM20}.

This paper is organized as follows. The
CDO-Fb scheme for the space discretization is presented in \cref{sec:spat_disc}.
The fully discrete schemes are introduced in \cref{sec:disc_pb} using either
the monolithic or the AC approaches and either an implicit or an explicit treatment
of the convective term. Numerical results are discussed in \cref{sec:num_res}
on two- and three-dimensional test cases. Finally, conclusions are drawn in \cref{sec:conc}.

\section{Space discretization by the CDO-Fb scheme}
\label{sec:spat_disc}

In this section, we recall the 
CDO-Fb scheme for the space discretization of the steady
incompressible Navier--Stokes equations introduced in \cite{BoEM20}.
In weak form, and considering homogeneous Dirichlet conditions for
simplicity (and fixing the pressure mean value to zero), one seeks
$(\Vex,\Pex)\in \Hilbertv1_0(\Dom)\times \Lebesgue2_0(\Dom)$ (standard notation for 
Lebesgue and Sobolev spaces is employed) such that
\begin{equation}
  \label{eq:NS-steady}
  \left\{
    \begin{alignedat}{2}
      \Visc a(\Vex,\Vtest) + t(\Vex;\Vex,\Vtest) + b(\Vtest,\Pex) ={}&l(\Vtest)\,,&\qquad &\forall \Vtest\in \Hilbertv1_0(\Dom)\,,\\
      b(\Vex,\Ptest)                                              ={}&0\,,&\qquad&\forall \Ptest\in \Lebesgue2_0(\Dom)\,,
    \end{alignedat}
  \right. 
\end{equation}
with the bilinear and trilinear forms 
\begin{equation} \label{eq:exact_forms}
  a(\Vex,\Vtest)\coloneqq \int_\Dom \tgrad\Vex : \tgrad\Vtest\,, \qquad
  b(\Vtest,\Pex)\coloneqq -\int_\Dom (\grad{\cdot}\Vtest)\Pex\,, \qquad
  t(\Vttest;\Vex,\Vtest)\coloneqq \int_\Dom \big( (\Vttest{\cdot}\grad)\Vex\big){\cdot}\Vtest\,,
\end{equation}
and the linear form $l(\Vtest)\coloneqq \int_\Dom \Foex{\cdot}\Vtest$.

\subsection{Mesh entities and degrees of freedom}

The mesh discretizing $ \Dom $ is a finite collection $ \cellh
\coloneqq \Set{\cell} $ of nonempty, disjoint, open, polytopal subsets of
$\Real^{\dmn}$, usually referred to as cells. The mesh faces are
assumed to be planar and are gathered in the set $ \faceh \coloneqq \Set{\face} $ which can be subdivided
in the two disjoint sets $ \faceh^b \coloneqq \Set{ \face\, |\, \face\subset\pOm} $
collecting the mesh boundary faces and $ \faceh^i \coloneqq \faceh\setminus\faceh^b $ 
collecting the mesh interfaces. 
One associates with each mesh face $ \face $ a normal vector $\norma{\face}$ such that
if $\finf^b$, $\norma{\face}$ points outward $ \Dom $ and, if $ \finf^i $, the
direction is chosen arbitrarily and fixed once and for all. 
For a generic mesh entity $\elem=\cell$ or
$\elem=\face$, $\pt\ue$ denotes its barycenter, $\meas{\elem}$ its measure and
$\diam\ue\coloneqq\max_{\pt_1,\pt_2\in\elem} \meas{\pt_1-\pt_2}$ its diameter.
Moreover, $ \diam \coloneqq \max_{\cinc} \diam\uc $ is called the size of the mesh. We
consider the shape-regularity setting of \cite{DPEr15} for polytopal mesh
families, and we additionally require that every mesh entity
$\elem=\cell$ or $\elem=\face$ is star-shaped with respect to its barycenter
$\pt\ue$.
The faces composing the boundary of a cell $\cinc$ are collected in the set
$\faceh\uc \coloneqq \{ \face\in\faceh \tq \face\subset\bnd\cell \}$. For each face
$\finc$, a unit normal vector pointing outward $\cell$ is denoted by $\nfc =
\pm\norma{\face} $, depending on the orientation chosen for
$\norma{\face}$. We are going to make use of a subdivision of the cell $ \cell $ as
$\Pvol\uc \coloneqq \{\pfc\}_{\finc} $, where the subsets $\pfc$ are the nonempty,
disjoint subpyramids (or subtriangles if $\dmn=2$) obtained by considering the cell
barycenter $\pt\uc$ as apex, and a face $\finc$ as basis.

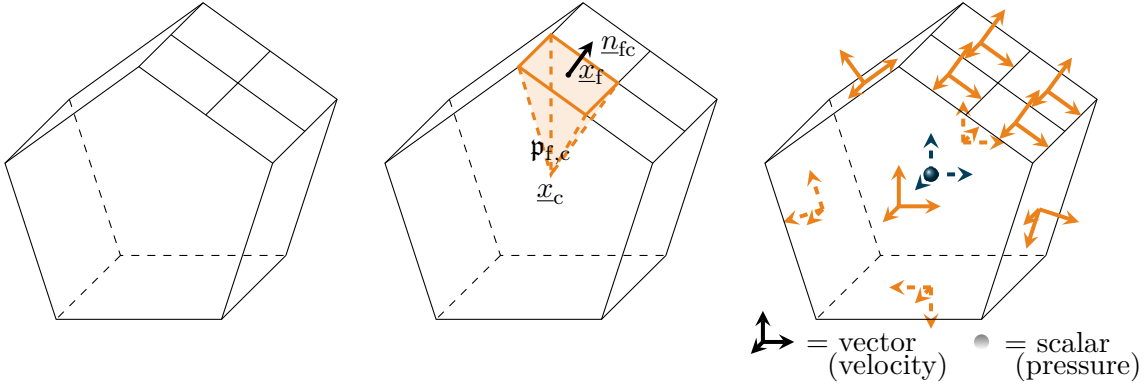
\begin{figure}
  \centering
  \begin{tikzpicture}[arr/.style={line width=1.5pt,->,>=stealth,color=#1},arr/.default={black},]
    \begin{scope}[xshift=0cm]
      \basepenta{} \pentainside{} \pentahang{}
    \end{scope}
    \begin{scope}[xshift=5cm]
      \basepenta{} \pentainside{} \pentahang{}
      \coordinate (HNS1) at ($(PH5)!.5!(P11)$);
      \draw[very thick,dashed,color=\pentacolor] (HN1) -- (C) -- (PH5);
      \draw[very thick,color=\pentacolor,fill=\pentacolor,draw opacity=1,fill opacity=0.15] (HN1) -- (P11) -- (PS5) -- (PH5) -- cycle;
      \draw[very thick,dashed,color=\pentacolor,fill=\pentacolor,draw opacity=1,fill opacity=0.15] (P11) -- (C) -- (PS5);
      \fill[\pentacolor,fill opacity=0.15] (PH5) -- (C) -- (PS5);
      \node[below] at (C) {$\pt\uc$};
      \node[above] at (C) {$\pvol_{\face,\cell}$};
      \pgfmathsetmacro{\ar}{0.55}
      \pgfmathparse{90-\intang/2}
      \pgfmathsetmacro\ca{\ar*cos(\pgfmathresult)}
      \pgfmathsetmacro\sa{\ar*sin(\pgfmathresult)}
      \draw[fill] (HNS1) circle (0.03) node[right] {$\pt\uf$};
      \draw[very thick,arr] (HNS1) -- node[pos=0.85,above,right] {$\nfc$} +(\ca,\sa,0);
    \end{scope}
    \begin{scope}[xshift=10cm]
      \basepenta{} \pentainside{} \pentahang{}
      \coordinate (HNS1) at ($(PH5)!.5!(P11)$);
      \coordinate (HNS2) at ($(PH5)!.5!(P15)$);
      \coordinate (HNS3) at ($(PH5)!.5!(P0.corner 1)$);
      \coordinate (HNS4) at ($(PH5)!.5!(P0.corner 5)$);
      \pgfmathsetmacro{\ar}{0.55}
      \foreach \v in {1,4} {
        \pgfmathparse{(90-\intang/2)+\intang*mod(\v,5)}
        \pgfmathsetmacro\ca{\ar*cos(\pgfmathresult)}
        \pgfmathsetmacro\sa{\ar*sin(\pgfmathresult)}
        \begin{pgfonlayer}{fg}
          \draw[arr=\pentacolor] (PH\v) -- +(\ca,\sa,0);%
          \draw[arr=\pentacolor] (PH\v) -- +(\sa,-\ca,0);%
          \draw[arr=\pentacolor] (PH\v) -- +(0,0,\ar);
        \end{pgfonlayer}%
      }
      \foreach \v in {2,3} {
        \pgfmathparse{(90-\intang/2)+\intang*mod(\v,5)}
        \pgfmathsetmacro\ca{\ar*cos(\pgfmathresult)}
        \pgfmathsetmacro\sa{\ar*sin(\pgfmathresult)}
        \begin{pgfonlayer}{bg}
          \draw[dashed,arr=\pentacolor] (PH\v) -- +(\ca,\sa,0);%
          \draw[dashed,arr=\pentacolor] (PH\v) -- +(\sa,-\ca,0);%
          \draw[dashed,arr=\pentacolor] (PH\v) -- +(0,0,\ar);
        \end{pgfonlayer}%
      }
      \pgfmathparse{(90-\intang/2)}
      \pgfmathsetmacro\ca{\ar*cos(\pgfmathresult)}
      \pgfmathsetmacro\sa{\ar*sin(\pgfmathresult)}
      \foreach \v in {1,...,4} {
         \begin{pgfonlayer}{fg}
          \draw[arr=\pentacolor] (HNS\v) -- +(\ca,\sa,0);%
          \draw[arr=\pentacolor] (HNS\v) -- +(\sa,-\ca,0);%
          \draw[arr=\pentacolor] (HNS\v) -- +(0,0,\ar);
        \end{pgfonlayer}
      }
      \begin{pgfonlayer}{bg}
        \draw[dashed,arr=\pentacolor] (O0) -- +(\ar,0,0);
        \draw[dashed,arr=\pentacolor] (O0) -- +(0,\ar,0);
        \draw[dashed,arr=\pentacolor] (O0) -- +(0,0,-\ar);
      \end{pgfonlayer}
      \begin{pgfonlayer}{fg}
        \draw[arr=\pentacolor] (O1) -- +(\ar,0,0);
        \draw[arr=\pentacolor] (O1) -- +(0,\ar,0);
        \draw[arr=\pentacolor] (O1) -- +(0,0,\ar);
      \end{pgfonlayer}
      \begin{pgfonlayer}{bg}
        \draw[dashed,arr=\cllcolor] (C) -- +(\ar,0,0);
        \draw[dashed,arr=\cllcolor] (C) -- +(0,\ar,0);
        \draw[dashed,arr=\cllcolor] (C) -- +(0,0,\ar);
        \shade[ball color=\presscolor] (C) circle (0.10);
      \end{pgfonlayer}

      \pgfmathsetmacro{\blw}{-0.3}
      \pgfmathsetmacro{\lft}{-0.7}
      \coordinate (BP13) at ($ (P13) + (\lft, \blw) $);
      \draw[arr] (BP13) -- +(\ar*.75,0,0) node[right] {= vector};
      \node[below right] at ($(BP13)+(\ar*.75,0,0)$) {\hphantom{=}(velocity)};
      \draw[arr] (BP13) -- +(0,\ar*.75,0);
      \draw[arr] (BP13) -- +(0,0,\ar*.75);

      \coordinate (BP14) at ($ (P14) + (0, \blw) $);
      \shade (BP14) circle (0.10) node[right] {\hspace*{1ex}{= scalar}};
      \node[below right] at (BP14) {\hphantom{=}(pressure)};
    \end{scope}
  \end{tikzpicture}
  \caption{Mesh, notation and DoFs. Left: Example of a polyhedral cell with a hanging
    node. Center: The barycenter $\pt\uc$ of the cell
    $\cell$ and a face $\face\in\faceh_\cell$ with its barycenter, $\pt\uf$, its normal
    $\nfc$ and its sub-pyramid $\pfc$. 
    Right: Full set of velocity and pressure DoFs for the considered cell (dashed arrows are used for velocity DoFs on hidden faces).}
  \label{fig:cell}
\end{figure}

For a generic mesh entity $\elem=\cell$ or $\elem=\face$, ${\Poly}^q(\elem)$ 
(resp., $\vect{\Poly}^q(\elem)$ and $\tens{\Poly}^q(\elem)$) 
is composed of the restriction to $\elem$
of the scalar-valued (resp., $\Real^{\dmn}$-valued
and $\Real^{\dmn\times\dmn}$-valued) polynomials of degree at most $ q $. 
Moreover, for a collection $\elemh$ of mesh entities,
$\Poly^{q}(\elemh)$, $\vect{\Poly}^{q}(\elemh)$, and $\tens{\Poly}^{q}(\elemh)$ 
refer to piecewise-polynomial functions; for instance, $
\Poly^1(\faceh\uc) \coloneqq \bigtimes_{\finc} \Poly^1(\face) $. Notice also that for
$\elem=\cell$ or $\elem=\face $, $\Poly^0(\elem)\equiv\Real$.

In the CDO-Fb framework, the velocity is hybrid, meaning that it has cell- and
face-based DoFs. Hence, the global velocity space is
\begin{equation}
\hvsp\uh\coloneqq \bigtimes_{\finf} \vect{\Poly}^0(\face) \times \bigtimes_{\cinc} \vect{\Poly}^0(\cell)\,.
\end{equation}
A generic element of $ \hvsp\uh $ is denoted by $ \hVtest\uh \coloneqq \big( ( \Vtest\uf 
)_{\finf}, ( \Vtest\uc )_{\cinc} \big) $.
Notice that the velocity DoFs at the interfaces are uniquely defined. 
Moreover, the velocity DoFs
associated with a generic cell $ \cell\in \cellh $ 
are denoted by 
\begin{equation} 
\hVtest\uc \coloneqq \big( (
\Vtest\uf )_{\finc}, \Vtest\uc \big) \in \hvsp\uc \coloneqq \bigtimes_{\finc}
\vect{\Poly}^0(\face) \times \vect{\Poly}^0(\cell)\,.
\end{equation} 
The pressure has only cell-based DoFs, so that the global pressure space is
\begin{equation} 
\spsp\uh \coloneqq \bigtimes_{\cinc} \spsp\uc \ni \Pex\uh \coloneqq \left( \Pex\uc \right)_{\cinc}\,,\qquad \spsp\uc \coloneqq \Poly^0(\cell)\,.
\end{equation}
In order to account for the
velocity boundary conditions and the constraint on the pressure average, we consider the
subspaces 
\begin{equation}
\hvspz\coloneqq \Set{\hVtest\uh\in\hvsp\uh\,|\,
\Vtest\uf=\vect{0}\ \forall\face\in\faceh^b}, \qquad
\spsps\coloneqq \Set{\Pex\uh\in\spsp\uh\,|\, \sum_{\cinc} \meas{\cell}\Pex\uc = 0}.
\end{equation}

Finally, for $\elem=\cell$ or $\elem=\face$, $\prj\ue$ (resp., $\vprj\ue$) denotes the
$\Lebesgue2$-orthogonal projection onto $\Poly^0(\elem)$ (resp., $\vect{\Poly}^0(\elem)$). 
The $\Lebesgue2$-orthogonal projection onto $\spsp\uh$ is defined cellwise
so that $\prj\uh(\Ptest)\coloneqq\left(\prj\uc(\Ptest\restr{\cell})\right)_{\cinc}$ for all 
$\Ptest\in\Lebesgue1(\Dom)$, whereas the projection onto the 
hybrid discrete space $\hvsp\uh$ is defined such that
$\hvprj\uh(\Vtest)\coloneqq\big((\vprj\uf(\Vtest\restr{\face}))_{\finf},
(\vprj\uc(\Vtest))_{\cinc}\big)$ for all $\vect{v}\in \vect{H}^{s}(\Dom)$, $s>\frac12$.
Similarly, for the local hybrid space $\hvsp\uc$, we set
$\hvprj\uc(\Vtest)\coloneqq\big((\vprj\uf(\Vtest\restr{\face}))_{\finc},
\vprj\uc(\Vtest)\big)$ for all $\vect{v}\in \vect{H}^{s}(\cell)$.

\subsection{Discrete diffusion-like bilinear form}

Consider a cell $\cinc$ and the velocity DoFs $\hvsp\uc$ 
associated with $\cell$. We define locally in the mesh cell $\cell$ a
velocity gradient reconstruction operator, $\tGh\uc:\hvsp\uc\to\tens{\Poly}^0(\Pvol\uc)$, i.e., for all $\hVtest\uc\in \hvsp\uc$, $\tGh\uc(\hVtest\uc)$ is 
a piecewise constant tensor-valued field on
the subpyramids collected in $\Pvol\uc$. In the context
of CDO schemes, this type of operator was introduced in \cite{Bone14}; see also 
\cite{BDPE15}. It is composed of a consistent part, $\tGhavg\uc$, plus a
stabilization part, $\tGhstb\uc$, so that
\begin{equation}
\tGh\uc(\hVtest\uc)\coloneqq \tGhavg\uc(\hVtest\uc)+\tGhstb\uc(\hVtest\uc)\,,
\end{equation}
where $\tGhavg\uc(\hVtest\uc)\in \tens{\Poly}^0(\cell)\subset \tens{\Poly}^0(\Pvol\uc)$ and $\tGhstb\uc(\hVtest\uc)\in \tens{\Poly}^0(\Pvol\uc)$ are defined as follows:
\begin{equation}
  \begin{aligned}
\tGhavg\uc(\hVtest\uc)\coloneqq{}&\frac1{\meas{\cell}}\sum_{\finc}\meas{\face}(\Vtest\uf-\Vtest\uc)\otimes\nfc\,,\\
\tGhstb\uc(\hVtest\uc)\restr{\pfc} \coloneqq{}&
    \CDOstb\frac{\meas{\face}}{\meas{\pfc}} \Big( (\Vtest\uf - \Vtest\uc) - \tGhavg\uc(\hVtest\uc)(\pt\uf - \pt\uc) \Big) \otimes \nfc\,, \qquad \forall \face\in \faceh\uc\,.
  \end{aligned}
\end{equation}
The positive stabilization parameter $\CDOstb$ is user-defined (with the only requirement of being positive): choosing $\CDOstb\coloneqq1$ recovers
the generalized Crouzeix--Raviart scheme from \cite{DPLe15}, whereas
$\CDOstb\coloneqq\frac1{\sqrt{\dmn}}$ leads to the HFV scheme from
\cite{EyGH10}. 

The global version of the gradient reconstruction operator, $\tGh\uh:\hvsp\uh
\to \bigtimes_{\cinc}\tens{\Poly}^0(\Pvol\uc)$, is
defined cellwise so that, for any $\hVtest\uh\in\hvsp\uh$ and any $\cinc$, we have
$\tGh\uh(\hVtest\uh)\restr{\cell}\coloneqq\tGh\uc(\hVtest\uc)$. With the above operators 
in hand, the discrete diffusion-like bilinear form $a_h:\hvsp\uh\times \hvsp\uh\to
\mathbb{R}$ is defined as follows:
\begin{equation}
a_h(\hVex\uh,\hVtest\uh) \coloneqq \int_\Dom \tGh\uh(\hVex\uh) : \tGh\uh(\hVtest\uh)
= \sum_{\cinc} \int_\cell \tGh\uc(\hVex\uc):\tGh\uc(\hVtest\uc)\,.
\end{equation}

The above-defined gradient reconstruction operators enjoy several important
properties. First, $\tGh\uc$ is consistent for affine functions, meaning that
\begin{equation}
\tGh\uc\big(\hvprj\uc(\vect{v})\big) = \tens{\nabla}\vect{v}\,,\qquad \forall 
\vect{v}\in\vect{\Poly}^1(\cell)\,.
\end{equation}
(More precisely, we have $\tGhavg\uc\big(\hvprj\uc(\vect{v})\big) = \tens{\nabla}\vect{v}$
and $\tGhstb\uc\big(\hvprj\uc(\vect{v})\big) = \tens{0}$.) Moreover, the consistent and
stabilization parts of the reconstructed gradient are $\tens{L}^2(\cell)$-orthogonal,
\begin{equation}
\int_\cell \tGhavg\uc(\hVtest\uc) : \tGhstb\uc(\hVtest\uc)=0\,,\qquad \forall
\hVtest\uc \in \hvsp\uc\,.
\end{equation}
Finally, and most importantly, the following stability and boundedness properties hold true:
There is $\delta>0$ such that for all $\cinc$ and all $\hVtest\uc \in \hvsp\uc$,
\begin{equation} \label{eq:stab_Gc}
\delta \|\hVtest\uc\|_{1,\cell}^2 \le \big\|\tGh\uc(\hVtest\uc)\big\|_{\tens{L}^2(\cell)}^2
\le \delta^{-1} \|\hVtest\uc\|_{1,\cell}^2\,,
\end{equation}
with the discrete $H^1$-like seminorm defined on $\hvsp\uc$ such that
$\|\hVtest\uc\|_{1,\cell}^2\coloneqq \sum_{\face\in\faceh\uc} h_\cell^{-1} |\face| \, 
|\Vtest\uf-\Vtest\uc|^2$. The lower bound~\eqref{eq:stab_Gc} is the main reason for
introducing the stabilization part of the reconstructed gradient. Notice 
that $\|\hVtest\uh\|_{1,\mathrm{h}}^2\coloneqq \sum_{\cinc} \|\hVtest\uc\|_{1,\cell}^2$
defines a norm on $\hvspz$.

\subsection{Discrete velocity-pressure coupling}

The discrete velocity-pressure coupling hinges on a discrete divergence operator
$\Dh\uh\colon\hvsp\uh\to\spsp\uh$ that is defined cellwise. For all $\cinc$, the
local discrete divergence operator $\Dh\uc\colon\hvsp\uc\to\spsp\uc$ is such that
for all $\hVtest\uc\in \hvsp\uc$,
\begin{equation}
  \Dh\uc(\hVtest\uc)\coloneqq\trace\left(\tGhavg\uc(\hVtest\uc)\right)
  =\frac1{\meas{\cell}}\sum_{\finc}\meas{\face}(\Vtest\uf-\Vtest\uc)\cdot\nfc\,.
\end{equation}
The global operator is then defined by setting
$\Dh\uh(\hVtest\uh)\restr{\cell}=\Dh\uc(\hVtest\uc)$
for all $\hVtest\uh \in \hvsp\uh$ and all $\cinc$. Notice that
$\Dh\uh(\hVtest\uh)\in \spsps$ for all $\Vtest\uh\in \hvspz$.

With the above divergence operator in hand, the discrete bilinear form 
$b_h:\hvsp\uh\times \spsp\uh\to
\mathbb{R}$ handling the velocity-pressure coupling is defined as follows:
\begin{equation}
b\uh(\hVtest\uh,\Ptest\uh) \coloneqq 
- \int_\Dom \Dh\uh(\hVtest\uh) \Ptest\uh = - \sum_{\cinc} 
\meas{\cell}\Dh\uc(\hVtest\uc)\Ptest\uc\,.
\end{equation}
The same discrete bilinear form 
is found in HMM schemes \cite{DrEF15,DrEy17} and the lowest-order HHO scheme 
\cite{ABDP15,DELS16}. A crucial property is the following discrete inf-sup
condition: There is $ \isprm > 0 $, only depending on the mesh shape-regularity,
such that
\begin{equation}
\label{eq:infsup_disc}
\inf_{\Ptest\uh\in\spsps} \sup_{\hVtest\uh\in\hvspz}
\frac{|b\uh(\hVtest\uh,\Ptest\uh)|}{\nrm{\Ptest\uh}\uh\nrm{\hVtest\uh}_{1,\dsc}}\geq\isprm\,,
\end{equation}
where $\nrm{\hVtest\uh}_{1,\dsc}$ is defined above whereas 
$\nrm{\Ptest\uh}\uh^2\coloneqq\nrm{\Ptest\uh}_{\Lebesgue2(\Dom)}^2 =
\sum_{\cinc}\meas{\cell}\abs{\Ptest\uc}^2$ for all $\Ptest\uh\in\spsp\uh$. 
Another important property of the discrete divergence operator is its commuting
property with the $L^2$-orthogonal projection in the sense that 
$\Dh\uc\big(\hvprj\uc(\Vtest)\big)=\prj\uc(\dive\Vtest)$ for all $\cinc$
and all $\Vtest\in\vect{H}^{1}(\cell)$.

\subsection{Discrete convection operator}

Let us finally devise a discrete CDO-Fb counterpart of the trilinear form 
$t$ defined in~\eqref{eq:exact_forms}. To this purpose, we define the 
discrete trilinear form $t_h:\hvsp\uh\times \hvsp\uh\times \hvsp\uh\to \mathbb{R}$
such that 
\begin{equation}
  \label{eq:conv_def_glob}
  \advTri\uh(\hVttest\uh;\hVex\uh,\hVtest\uh)\coloneqq
  \sum_{\cinc} \advTri\uc(\hVttest\uc;\hVex\uc,\hVtest\uc)
  +\sum_{\finf^b}\advTri\uf(\Vttest\uf;\Vex\uf,\Vtest\uf)\,,
\end{equation}
with
\begin{equation}
  \label{eq:conv_def_other}
  \begin{aligned}
    \advTri\uc(\hVttest\uc;\hVex\uc,\hVtest\uc) \coloneqq
    {}&\frac{1}{2}\sum_{\finc}\meas{\face} (\Vttest\uf\cdot\nfc) (\Vex\uf-\Vex\uc)\cdot(\Vtest\uf+\Vtest\uc)\,,\\
    \advTri\uf(\Vttest\uf;\Vex\uf,\Vtest\uf) \coloneqq
    {}&\meas{\face} (\Vttest\uf\cdot\nfc)^- \Vex\uf\cdot\Vtest\uf\,,
  \end{aligned}
\end{equation}
where $(x)^-\coloneqq \frac12(|x|-x)$ denotes the negative part of any real number 
$x\in\mathbb{R}$. The trilinear form $\advTri\uf$ is related to the weak enforcement of
inflow boundary conditions and vanishes whenever at least one of the arguments
is in $\hvspz$. The discrete
trilinear form $t_h$ is inspired from the lowest-order HHO method
for scalar advection-diffusion equations introduced in \cite{DPED15}; see also
\cite[Remark 9]{BDPD19} for the Navier--Stokes equations. 

The first important property of the discrete trilinear form $t_h$ is its 
positivity and skew-symmetry. Indeed, assume that $\hVttest\uh\in\hvsp\uh$ is 
discretely divergence-free, i.e., $\Dh\uc(\hVttest\uc)=0$ for all $\cinc$.
Then, one readily verifies that
\begin{equation}
  \label{eq:conv_pos}
  \advTri\uh(\hVttest\uh;\hVex\uh,\hVex\uh) \geq 0\,, \qquad \forall \hVex\uh\in\hvsp\uh\,.
\end{equation}
Moreover, if additionally the normal component of $\hVttest\uh$ vanishes at the
boundary, namely $\Vttest\uf\cdot\norma{\face}=0$ for all $\finf^b$, then 
$\advTri\uh(\hVttest\uh;\cdot,\cdot)$ is skew-symmetric, i.e.
\begin{equation}
  \label{eq:conv_skew}
  \advTri\uh(\hVttest\uh;\hVex\uh,\hVex\uh)=0\,, \qquad\forall\hVex\uh\in\hvsp\uh\,.
\end{equation}
Property \eqref{eq:conv_pos} is crucial in the context of the Navier--Stokes
equations since it is instrumental to establish  the dissipativity of the
discrete problem. Notice also that \eqref{eq:conv_skew} is reminiscent of the
so-called Temam's trick on the discrete trilinear form. 
A further relevant property of the discrete trilinear form $t_h$ is its
limit-conformity.
This notion is described, e.g., in \cite{EyFG18}, and is verified for the above
CDO-Fb setting in \cite[Remark~2.50]{Mila20}; see also \cite[Prop.~6]{DPKr18}.

\begin{remark}[Upwinding]
An upwinding stabilization can be considered by adding to \eqref{eq:conv_def_glob}
the discrete trilinear form 
$\advTri\uh^{\rm u}(\hVttest\uh;\hVex\uh,\hVtest\uh) \coloneqq
\frac12\sum_{\finf^i}\sum_{\cinc\uf}\meas{\face} \abs{\Vttest\uf\cdot\nfc}(\Vex\uf-\Vex\uc)\cdot(\Vtest\uf-\Vtest\uc)$. This option is not considered further in the present work.
\end{remark}

\section{Time discretization by monolithic and artificial compressibility schemes}
\label{sec:disc_pb}

We present in this section the fully discrete problems obtained by using either a monolithic 
or an AC scheme for the time discretization and by considering either an
implicit or an explicit treatment of the convection operator. Moreover, we consider 
both first-order and second-order time schemes. 

For a hybrid velocity field $\hVtest\uh\coloneqq \big( ( \Vtest\uf 
)_{\finf}, ( \Vtest\uc )_{\cinc} \big)\in \hvsp\uh$, we denote its
cell-based components by $\Vtest\uC\coloneqq ( \Vtest\uc )_{\cinc}$.
Then, to discretize the time-derivative of the velocity, we only use
the cell-based components and consider the mass bilinear form such that 
\begin{equation}
m(\Vex\uC,\Vtest\uC)\coloneqq \sum_{\cinc} \meas{\cell} \Vex\uc{\cdot}\Vtest\uc\,.
\end{equation}
We also use the cell-based components to discretize the source term and we set
$l^n(\Vtest\uC) \coloneqq \sum_{\cinc} \int_{\cell} \Foex^n{\cdot}\Vtest\uc$ for all $n\ge1$.
For any hybrid velocity field
$\hVtest\uh\in\hvspz$, we define its discrete kinetic energy as
\begin{equation}
\label{eq:kin_en}
\EnrKin(\hVtest\uh)\coloneqq\frac12\nrm{\Vtest\uC}_{\Lebesguev2(\Dom)}^2 = \frac12 \sum_{\cinc}\meas{\cell}\abs{\Vtest\uc}_2^2\,.
\end{equation}

\subsection{First-order schemes}

\subsubsection{Monolithic scheme}

Let us first consider an implicit treatment of the convection term
which is dealt with by means of a Picard algorithm. Then, the 
first-order monolithic time-stepping scheme reads as follows.
For all $ n=1,\ldots,N $, iterate on $ k\geq1 $ until convergence:
Find $(\hVex\uh^{n,k},\Pex\uh^{n,k})\in\hvspz\times\spsps$ such that
\begin{equation}
  \label{eq:mono_fst_impl}
  \left\{
    \begin{aligned}
      \frac{1}{\deltat}m(\Vex\uC^{n,k}-\Vex\uC\nm[1],\Vtest\uC) + \Visc a\uh(\hVex\uh^{n,k},\hVtest\uh)
      + \advTri\uh(\hVex\uh^{n,k-1};\hVex\uh^{n,k},\hVtest\uh)
      +b\uh(\hVtest\uh,\Pex\uh^{n,k}) ={}& l^n(\Vtest\uC)\,,\\
      b\uh(\hVex\uh^{n,k},\Ptest\uh) ={}&\, 0\,.
    \end{aligned}
  \right.
\end{equation}
for all $\hVtest\uh\in\hvspz$ and all $\Ptest\uh\in\spsps$. Here, $\hVex\uh\nm[1]$
denotes the solution given by the Picard algorithm at the time step $ n-1 $. The
time-stepping is initialized with the initial condition
$\hVex\uh^0\coloneqq\hvprj\uh(\Vex_0)$. Moreover, at each time step, the Picard
iterations have to be initialized: a suitable choice is to take the solution at the
previous time step, i.e., $\hVex\uh^{n,0}\coloneqq\hVex\uh\nm[1]$.
Taking $\hVtest\uh=\hVex\uh^{n,k}$ in \eqref{eq:mono_fst_impl} and using simple arithmetic
manipulations together with $\advTri\uh(\hVex\uh^{n,k-1};\hVex\uh^{n,k},\hVex\uh^{n,k})=0$ 
owing to \eqref{eq:conv_skew}, we obtain the following discrete kinetic energy balance:
\begin{equation}
\label{eq:mono_kin_bal}
\EnrKin(\hVex\uh^{n,k})-\EnrKin(\hVex\uh\nm[1])+\EnrKin(\hVex\uh^{n,k}-\hVex\uh\nm[1])+\Visc\deltat\big\|{\tGh\uh(\hVex\uh^{n,k})}\big\|^2_{\Lebesguet{2}(\Dom)}=\deltat \, l^n(\Vex^{n,k}\uC)\,.
\end{equation}
Notice that the third and fourth terms on the left-hand side are non-negative, i.e.,
they only contribute to energy dissipation. 
Notice also that \eqref{eq:mono_kin_bal} remains valid as $k\to\infty$.

With an explicit treatment of the convection term, the 
first-order monolithic time-stepping scheme reads as follows:
For all $ n=1,\ldots,N $, find $(\Vex\uh^n,\Pex\uh^n)\in\hvspz\times\spsps$ such that
\begin{equation}
  \label{eq:mono_fst_expl}
  \left\{
    \begin{aligned}
      \frac{1}{\deltat}m(\Vex\uC^{n}-\Vex\uC\nm,\Vtest\uC) + \Visc a\uh(\hVex\uh^{n},\hVtest\uh)
      +b\uh(\hVtest\uh,\Pex\uh^{n})
      ={}&l^n(\Vtest\uC)- \advTri\uh(\hVex\uh\nm;\hVex\uh\nm,\hVtest\uh)\,,\\
      b\uh(\hVex\uh^{n},\Ptest\uh) ={}& 0\,,
    \end{aligned}
  \right.
\end{equation}
for all $\hVtest\uh\in\hvspz$ and all $\Ptest\uh\in\spsps$. 
Unfortunately, a dissipative discrete kinetic energy balance cannot be derived for 
\eqref{eq:mono_fst_expl} because the term 
$\advTri\uh(\hVex\uh\nm;\hVex\uh\nm,\hVex\uh^n)$ cannot be given an a priori sign.
Thus, the scheme \eqref{eq:mono_fst_expl} is subject to a first-order CFL restriction
on the time step for its stability. This condition will be investigated numerically
in \cref{sec:num_res}.  

\subsubsection{AC scheme}

Let us first consider an implicit treatment of the convection term. 
As in the monolithic approach, we consider Picard iterations to deal with
the resulting nonlinearity. Then, the first-order AC time-stepping scheme reads as follows.
For all $ n=1,\ldots,N $, iterate on $ k\geq1 $ until convergence:
Find $\Vex\uh^{n,k}\in\hvspz$ such that for all $\hVtest\uh\in\hvspz$,
\begin{multline}
  \label{eq:ac_fst_impl}
  \frac{1}{\deltat}m(\Vex\uC^{n,k}-\Vex\uC\nm[1],\Vtest\uC) +
  \Visc\big(a\uh(\hVex\uh^{n,k},\hVtest\uh)
  +\ACprm d\uh(\hVex\uh^{n,k},\hVtest\uh)\big)+ \advTri\uh(\hVex\uh^{n,k-1};\hVex\uh^{n,k},\hVtest\uh)\\
  = l^n(\Vtest\uC)-b\uh(\hVtest\uh,\Pex\uh^{n,k-1}) \,,
\end{multline}
and then set
\begin{equation} \label{eq:ac_fst_impl_p}
  \Pex\uh^{n,k}=\Pex\uh^{n,k-1}-\Visc\ACprm\Dh\uh(\hVex\uh^{n,k})\,,
\end{equation}
recalling that $\eta>0$ is a user-defined parameter (choices are discussed 
in \cref{sec:num_res}). In~\eqref{eq:ac_fst_impl}, we introduced 
the discrete div-div bilinear form such that
\begin{equation}
d\uh(\hVex\uh,\hVtest\uh) \coloneqq \int_\Dom \Dh\uh(\hVex\uh)\Dh\uh(\hVtest\uh)
= \sum_{\cinc} \meas{c} \Dh\uc(\hVex\uc)\Dh\uc(\hVtest\uc)\,.
\end{equation}
Notice that in~\eqref{eq:ac_fst_impl}-\eqref{eq:ac_fst_impl_p},
we took advantage of the Picard iterations to update progressively the
pressure. Letting $k\to\infty$, we observe that the velocity and pressure
updates stemming from the monolithic and the AC schemes are the same.
Differently from the monolithic coupling, the initialization of the Picard
iteration requires to specify the velocity and the pressure, which are here
taken from the previous time step or the initial condition.

With an explicit treatment of the convection term, the first-order AC scheme reads as
follows: For all $ n=1,\ldots,N $, find $(\Vex\uh^n,\Pex\uh^n)\in\hvspz\times\spsps$ 
such that for all $\hVtest\uh\in\hvspz$,
\begin{multline}
  \label{eq:ac_fst_expl}
  \frac{1}{\deltat}m(\Vex\uC^{n}-\Vex\uC^{n-1},\Vtest\uC) + \Visc\big(a\uh(\hVex\uh^{n},\hVtest\uh)
  +\ACprm d\uh(\hVex\uh^{n},\hVtest\uh)\big)\\
  = l^n(\Vtest\uC)-b\uh(\hVtest\uh,\Pex\uh^{n-1})
  -\advTri\uh(\hVex\uh^{n-1};\hVex\uh^{n-1},\hVtest\uh)\,,
\end{multline}
and then set
\begin{equation}
  \Pex\uh^n=\Pex\uh\nm-\Visc\ACprm\Dh\uh(\hVex\uh^n)\,.
\end{equation}
As for the monolithic scheme, the explicit AC scheme is subject to a first-order
CFL restriction on the time step for its stability. This condition will be
investigated numerically in \cref{sec:num_res}.  


\begin{remark}[Initial pressure]
As is classical with the AC technique, an initial pressure needs to be specified.
In all the numerical experiments reported in this work, the initial pressure
is indeed known. If this were not the case, a possibility proposed in \cite{GuMi15} (see also
\cite[Remarks 74.4 and 75.7]{ErGu21III}) is to recover the
initial pressure by solving the following steady problem:
  \begin{equation}
    \left\{
      \begin{aligned}
        &\Delta\Pex^0={}\dive\Foex^0\,,\\
        &\pd{\Pex^0}{\norm}\restr{\pOm}=\left(\Foex^0-(-\Visc\lapl\Vex^0+(\Vex^0\cdot\grad)\Vex^0)\right)\restr{\pOm}\cdot\norma{\pOm}\,.
      \end{aligned}
    \right.
  \end{equation}
Obviously, the initial pressure is zero is the fluid is initially at rest and the source term is divergence-free and has a zero normal component on the boundary.
\end{remark}

\subsection{Second-order schemes}
\label{sec:second_order}

\subsubsection{Monolithic scheme}
One standard technique is to use a second-order backward differentiation formula
(BDF2) to devise second-order time-stepping schemes within the monolithic
approach. In general, BDF2 is employed for every time step $n\ge2$, and 
an implicit Euler step can still be considered for $n=1$.
With an implicit treatment of the convection term by means of a Picard algorithm,
the second-order monolithic scheme reads as follows. 
For all $ n=2,\ldots,N $, iterate on $ k\geq1 $ until convergence:
Find $(\hVex\uh^{n,k},\Pex\uh^{n,k})\in\hvspz\times\spsps$ such that
\begin{equation}
  \left\{
    \begin{aligned}
      &
      \begin{multlined}
        \frac{1}{2\deltat}m(3\Vex\uC^{n,k}-4\Vex\uC\nm[1]+\Vex\uC\nm[2],\Vtest\uC)
        +\Visc a\uh(\hVex\uh^{n,k},\hVtest\uh)\\+\advTri\uh(\hVex\uh^{n,k-1};\hVex\uh^{n,k},\hVtest\uh)
        +b\uh(\hVtest\uh,\Pex\uh^{n,k}) ={} l^n(\Vtest\uC)\,,
      \end{multlined}
      \\
      &b\uh(\hVex\uh^{n,k},\Ptest\uh) ={}\, 0\,,
    \end{aligned}
  \right.
\end{equation}
for all $\hVtest\uh\in\hvspz$ and all $\Ptest\uh\in\spsps$. Recall that
$\hVex\uh\nm[1]$ (resp.~$\hVex\uh\nm[2]$) denotes the solution
given by the Picard algorithm at the time step $n-1$ (resp.~$n-2$). The
initialization of the iterative procedure is done with an approximation of the
solution, for instance $\hVex\uh\nm$ (as for the first-order scheme), or by using the
second-order extrapolation formula $(2\hVex\uh\nm-\hVex\uh\nm[2])$.
Finally, for the BDF2 time discretization, it is possible to derive by means of
classical algebraic manipulations (detailed, e.g., in \cite[Lemma~68.1]{ErGu21III} 
for the heat equation) a time-discrete kinetic energy balance with dissipation, 
in the same spirit as in \eqref{eq:mono_kin_bal}.

With an explicit treatment of the convection term, the second-order 
monolithic scheme reads as follows:
For all $ n=2,\ldots,N $, find $(\hVex\uh^n,\Pex\uh^n)\in\hvspz\times\spsps$ solving
\begin{equation}
  \left\{
    \begin{aligned}
      &\begin{multlined}
        \frac{1}{2\deltat}m(3\Vex\uC^n-4\Vex\uC\nm+\Vex\uC\nm[2],\Vtest\uC) +
        \Visc a\uh(\hVex\uh^{n},\hVtest\uh) +b\uh(\hVtest\uh,\Pex\uh^{n})\\
        ={}l^n(\Vtest\uC)
        -\big(2\advTri\uh(\hVex\uh\nm;\hVex\uh\nm,\hVtest\uh)-\advTri\uh(\hVex\uh\nm[2];\hVex\uh\nm[2],\hVtest\uh)\big)\,,
      \end{multlined}\\
      &b\uh(\hVex\uh^{n},\Ptest\uh) ={} 0\,,
    \end{aligned}
  \right.
\end{equation}
for all $\hVtest\uh\in\hvspz$ and all $\Ptest\uh\in\spsps$.
Notice that a second-order extrapolation formula is used for the convective term on the
right-hand side; here, we applied the extrapolation formula to the operator, but 
it is also possible to consider applying the extrapolation formula to the discrete velocity
and then form the convection operator. Finally, as for the first-order time discretization,
a dissipative kinetic energy balance is not available with an explicit treatment of the convection
term, and a CFL restriction on the time step is required for stability.

\subsubsection{AC scheme}

It has been shown in \cite{GuMi15} that arbitrary order in time
can be achieved by combining the AC method with a bootstrap
technique. In order to obtain the $k$-th order, $k$ linear systems similar to
\eqref{eq:ust_ac} have to be solved. An alternative technique
to reach arbitrary order in time, which is also discussed in \cite{GuMi15}, is
based on a defect correction procedure, leading to a similar computational cost.
To avoid the proliferation of variants,
we focus only on the bootstrap technique. 

As motivated in the introduction, we only consider an explicit treatment of the
convection term. The second-order AC time-stepping scheme reads as follows:
For all $ n \geq 1 $, find $(\hVex_{1,\dsc}^n,\Pex_{1,\dsc}^n)\in\hvspz\times\spsps$ such that
for all $\hVtest\uh\in\hvspz$,
\begin{subequations}
  \label{eq:boot_expl}
  \begin{empheq}[left=\empheqlbrace]{align}
    &\begin{multlined}
      \frac{1}{\deltat}m(\Vex_{1,\cellh}^{n}-\Vex_{1,\cellh}^{n-1},\Vtest\uC)
      + \Visc \big(a\uh(\hVex_{1,\dsc}^{n},\hVtest\uh)
      + \ACprm d\uh(\hVex_{1,\dsc}^{n},\hVtest\uh)\big)\\
      = l^n(\Vtest\uC)-b\uh(\hVtest\uh,\Pex_{1,\dsc}^{n-1})
      -\advTri\uh(\hVex_{1,\dsc}^{n-1};\hVex_{1,\dsc}^{n-1},\hVtest\uh)\,,
    \end{multlined}\\
    &\Pex_{1,\dsc}^n = \Pex_{1,\dsc}\nm - \Visc\ACprm \Dh\uh(\hVex_{1,\dsc}^n),\quad \dACPex^n_{1,\dsc} \coloneqq \Pex_{1,\dsc}^n - \Pex_{1,\dsc}\nm \, .
  \end{empheq}
  Then, for all $n\geq2$, find
  $(\hVex_{2,\dsc}^n,\Pex_{2,\dsc}^n)\in\hvspz\times\spsps$ such that for all $\hVtest\uh\in\hvspz$,
  \begin{empheq}[left=\empheqlbrace]{align}
    & \begin{multlined}
      \frac{1}{2\deltat}m(3\Vex_{2,\cellh}^n-4\Vex_{2,\cellh}\nm+\Vex_{2,\cellh}\nm[2],\Vtest\uC)
      +\Visc\big(a\uh(\hVex_{2,\dsc}^{n},\hVtest\uh)
      + \ACprm d\uh(\hVex_{2,\dsc}^{n},\hVtest\uh)\big)\\
      = l^n(\Vtest\uC)-b\uh(\hVtest\uh,\Pex_{2,\dsc}^{n-1}+\dACPex^n_{1,\dsc})
      -\big(2\advTri\uh(\hVex_{2,\dsc}\nm;\hVex_{2,\dsc}\nm,\hVtest\uh)-\advTri\uh(\hVex_{2,\dsc}\nm[2];\hVex_{2,\dsc}\nm[2],\hVtest\uh)\big)\,,
    \end{multlined}\\
    &\Pex_{2,\dsc}^n = \Pex_{2,\dsc}\nm + \dACPex^n_{1,\dsc}  
    - \Visc\ACprm \Dh\uh(\hVex_{2,\dsc}^n)\,,
  \end{empheq}
\end{subequations}
with the following initialization choice: $\hVex_{2,\dsc}^1 \coloneqq
\hVex_{1,\dsc}^1$, $\hVex_{2,\dsc}^0 \coloneqq \hVex_{1,\dsc}^0 \coloneqq
\hVex_{\dsc}^0$ for the velocity and $\Pex_{2,\dsc}^1 \coloneqq
\Pex_{1,\dsc}^1$ for the pressure.

\section{Numerical results}
\label{sec:num_res}

In this section, we present numerical results to assess the performances of the 
AC method together with the CDO-Fb discretization. We first consider the 2D
Taylor--Green vortex test case  with the goal to verify the convergence rates
in time for the first- and second-order time-stepping schemes and both treatments of
the convection term. We also use this test case to study the CFL restriction on
the time step in the case of an explicit treatment of the convection term.
Then, we consider a modified 3D Taylor--Green vortex to compare all of the above
strategies in terms of efficiency, that is, we compare the reached error
thresholds on the velocity and the pressure with the needed CPU time. 
All the results are computed with the CDO module available in the open-source,
single-phase CFD solver \CS{} \cite{CS}.
For the 2D test cases, the linear systems are solved with a sparse
direct solver available in the MUMPS \cite{MUMPS} library. For the 3D test cases,
iterative solvers are considered and further discussed in Section~\ref{sec:3D}.

\subsection{2D Taylor--Green vortex}
The Taylor--Green vortex \cite{TaGr37} is a well-known 2D test case usually
considered to evaluate the performances of an unsteady Navier--Stokes solver.
In fact, it constitutes an analytic solution of the 2D Navier--Stokes given by
\begin{equation}
  \label{eq:TGV}
  \left\{
  \begin{aligned}
    \Vex_{\mathrm{TGV}}(x,y) \coloneqq{}& \exp(-2\Visc t)
    \begin{bmatrix} \sin(x)\cos(y) \\ -\cos(x)\sin(y) \end{bmatrix}\,,\\
    \Pex_{\mathrm{TGV}}(x,y) \coloneqq{}& \frac14\exp(-4\Visc t)(\cos(2x)+\cos(2y))\,.
  \end{aligned}
  \right.
\end{equation}
The computational domain is $ \Dom\coloneqq[0,\,2\pi]^2 $, whereas the observation time 
$T$ depends on the numerical experiment. The viscosity $\Visc$ is uniform but different values are
considered in order to modify the Reynolds number $\Rey$. Since the reference length and velocity
for this test case can be set to $ L\coloneqq1 $ and $U\coloneqq1$, we obtain 
$ \Rey=\frac{1}{\Visc} $.
Cartesian and polygonal meshes are considered. 

\subsubsection{Convergence results}
The convergence in time is first verified for all the schemes presented in
\cref{sec:disc_pb}. Let $\Vtest(t,\pt)$ (resp.~$\Ptest(t,\pt)$) be the analytical velocity
(resp.~pressure), and $\Set{\hVtest\uh^n}_{n=1,\ldots,N}\in[\hvspz]^{N}$
(resp.~$\Set{\Ptest\uh^n}_{n=1,\ldots,N}\in[\spsps]^{N}$) its space-time
discretization. The following space-time error norms are used: 
\begin{equation} \label{eq:def_errors_num}
  \begin{aligned}
    \nrm{\hvprj\uh(\Vtest)-\Set{\hVtest\uh^n}_n}_{\ell^2(\Lebesgue2)}^2\coloneqq&{}\sum_{n=1}^N\deltat\nrm{\hvprj\uh(\Vtest(t^n,\cdot))-\hVtest\uh^n}\uC^2\coloneqq\sum_{n=1}^N\deltat\sum_{\cinc}\meas{\cell}\abs{\vprj\uc(\Vtest(t^n,\cdot))-\hVtest\uc^n}^2\;,
    \\
    \nrm{\hvprj\uh(\Vtest)-\Set{\hVtest\uh^n}_n}_{\ell^2(\Hilbert1)}^2\coloneqq&{}\sum_{n=1}^N\deltat\sum_{\cinc}\sum_{\finc}\meas{\pfc}\abs{\tGh\uc\big(\hvprj\uc(\Vtest(t^n,\cdot))-\hVtest\uc^n\big)\restr{\pfc}}^2\;,
    \\
    \nrm{\vprj\uh(\Ptest)-\Set{\Ptest\uh^n}_n}_{\ell^2(\Lebesgue2)}^2\coloneqq&{}\sum_{n=1}^N\deltat\nrm{\prj\uh(\Ptest(t^n,\cdot))-\Ptest\uh^n}\uC^2\coloneqq\sum_{n=1}^N\deltat\sum_{\cinc}\meas{\cell}(\vprj\uc(\Ptest(t^n,\cdot))-\Ptest\uc^n)^2 \;.
  \end{aligned}
\end{equation}
These errors are then normalized by the value obtained for the corresponding norm
of the exact solution.

\begin{figure}
  \centering
  \begin{tikzpicture}
    \begin{groupplot}[
        group style={group size=3 by 2,
          horizontal sep=1cm,vertical sep=2mm,
          xlabels at=edge bottom,
        },
        width=0.35\linewidth,
        xmode=log,ymode=log,
        base plot,xlabel={$\deltat$},
        xtick=data,
      ]
      \nextgroupplot[table/x=dt,xticklabels={},table/y expr={\thisrow{Vel_l2tS_errcM}/sqrt(4.89420)},
        title={$\ell^2(\Lebesgue2)$ velocity error}]
        \addplot[mark=none,tgvMONO]   table {\MONOexp};\label{mrk:tgvMONO}
        \addplot[mark=none,tgvAC]   table {\MONOexp};\label{mrk:tgvAC}
        \addplot[expl,tgvMONO] table {\MONOexp};\label{mrk:tgvexpl}
        \addplot[expl,tgvAC]   table {\ACexp};
        \addplot[pica,tgvMONO] table {\MONOpica};\label{mrk:tgvpica}
        \addplot[pica,tgvAC]   table {\ACpica};
        \addplot[order]  expression[domain=0.01875:0.6] {0.06*x} node[ordernode,below] {1\fst{} ord};
      \nextgroupplot[table/x=dt,xticklabels={},table/y expr={\thisrow{Vel_h1tS_err_M}/sqrt(19.57676)},
        title={$\ell^2(\Hilbert1)$ velocity error}]
        \addplot[expl,tgvMONO] table {\MONOexp};
        \addplot[expl,tgvAC]   table {\ACexp};
        \addplot[pica,tgvMONO] table {\MONOpica};
        \addplot[pica,tgvAC]   table {\ACpica};
        \addplot[order]  expression[domain=0.01875:0.6] {0.2*x} node[ordernode] {1\fst{} ord};
      \nextgroupplot[table/x=dt,xticklabels={},table/y expr={\thisrow{Prs_l2tS_prsc}/sqrt(0.308404)},
        title={$\ell^2(\Lebesgue2)$ pressure error}]
        \addplot[expl,tgvMONO] table {\MONOexp};
        \addplot[expl,tgvAC]   table {\ACexp};
        \addplot[pica,tgvMONO] table {\MONOpica};
        \addplot[pica,tgvAC]   table {\ACpica};
        \addplot[order]  expression[domain=0.01875:0.6] {0.5*x} node[ordernode,below] {1\fst{} ord};
      %
      \nextgroupplot[table/x=dt,xticklabels={},table/y expr={\thisrow{Vel_l2tS_errcM}/sqrt(4.89420)},
        xticklabels={$\frac{T}{2}$,$\frac{T}{4}$,$\frac{T}{8}$,$\frac{T}{16}$,$\frac{T}{32}$,$\frac{T}{64}$},]
        \addplot[expl,tgvMONO] table {\BDFexp};
        \addplot[expl,tgvAC]   table {\BOOTexp};
        \addplot[pica,tgvMONO] table {\BDFpica};
        \addplot[order]  expression[domain=0.01875:0.6] {0.06*x*x} node[ordernode,below] {2\scn{} ord};
      \nextgroupplot[table/x=dt,xticklabels={},table/y expr={\thisrow{Vel_h1tS_err_M}/sqrt(19.57676)},
        xticklabels={$\frac{T}{2}$,$\frac{T}{4}$,$\frac{T}{8}$,$\frac{T}{16}$,$\frac{T}{32}$,$\frac{T}{64}$},]
        \addplot[expl,tgvMONO] table {\BDFexp};
        \addplot[expl,tgvAC]   table {\BOOTexp};
        \addplot[pica,tgvMONO] table {\BDFpica};
        \addplot[order]  expression[domain=0.01875:0.6] {0.2*x*x} node[ordernode] {2\scn{} ord};
      \nextgroupplot[table/x=dt,xticklabels={},table/y expr={\thisrow{Prs_l2tS_prsc}/sqrt(0.308404)},
        xticklabels={$\frac{T}{2}$,$\frac{T}{4}$,$\frac{T}{8}$,$\frac{T}{16}$,$\frac{T}{32}$,$\frac{T}{64}$},]
        \addplot[expl,tgvMONO] table {\BDFexp};
        \addplot[expl,tgvAC]   table {\BOOTexp};
        \addplot[pica,tgvMONO] table {\BDFpica};
        \addplot[order]  expression[domain=0.01875:0.6] {0.5*x*x} node[ordernode,below] {2\scn{} ord};
    \end{groupplot}
  \end{tikzpicture}
  \begin{tabular}{r@{:~}ll}
    Coupling & Monolithic\hspace{1em}\ref{mrk:tgvMONO} & AC($\ACprm=10$)\hspace{1em}\ref{mrk:tgvAC}\\
    Convection & Implicit \ref{mrk:tgvpica}
               & Explicit \ref{mrk:tgvexpl}
  \end{tabular}
  \caption{2D Taylor--Green vortex. Convergence results for the space-time
    velocity and pressure errors. $\Rey=1$, $T=1.2$. Top: first-order
    schemes; bottom: second-order schemes. Cartesian mesh composed of
  $512^2$ cells.}
  \label{fig:tgv_cvg}
\end{figure}
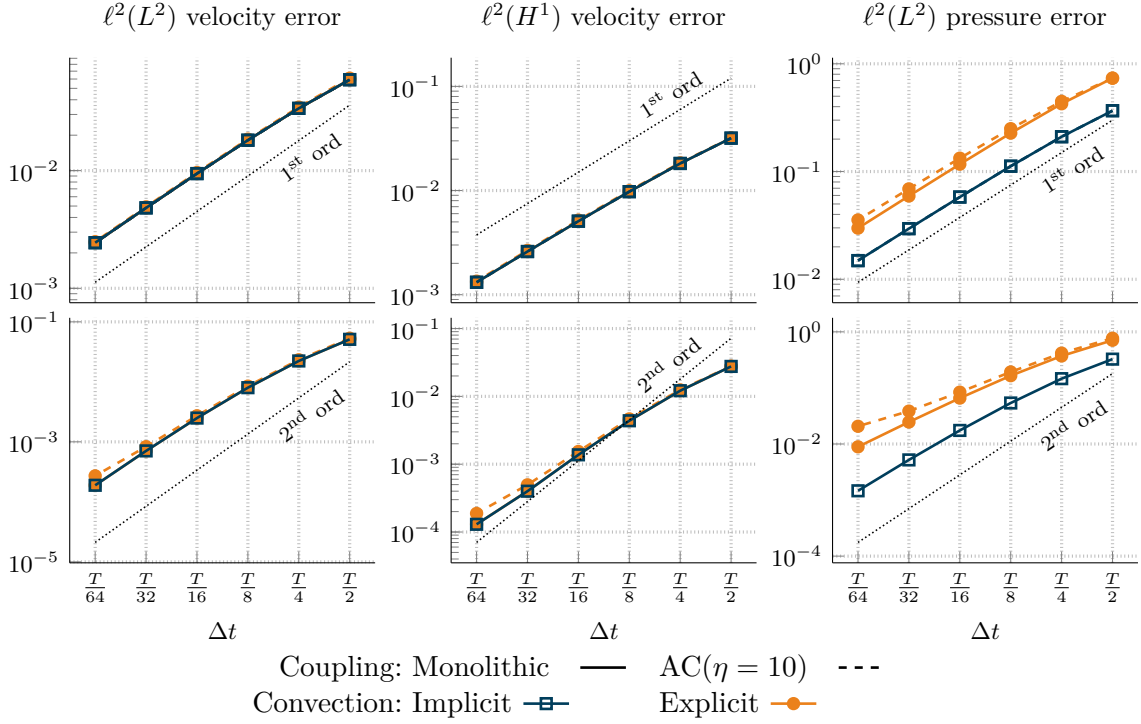

\Cref{fig:tgv_cvg} shows the results for first- (top row) and second-order
time-stepping schemes (bottom row) on a Cartesian mesh composed of $512^2$ cells. The Reynolds
number is first set to $\Rey=1$ (i.e., we take $\nu=1$), 
and the AC parameter is set to $\ACprm=10$. We consider both the
implicit and the explicit treatments of the convection term, and observe that in the
present setting, no instability is observed for the considered values of the time
step. This favorable situation is attributed to the moderate value of the Reynolds
number, and we will see below (see \cref{ssub:tgv_stab}) that a CFL stability
restriction needs to be enforced for larger Reynolds numbers. The results presented
in \Cref{fig:tgv_cvg} show that the monolithic and AC approaches lead to very
similar errors, both for the velocity and the pressure. The convection treatment has
an effect only on the pressure error (right column): the implicit treatment (squares)
is more accurate than the explicit one (circles). Indeed, due to the peculiar
construction of the Taylor--Green vortex for which the pressure gradient and the
convection term balance out, a good approximation of the convection term leads also
to an accurate approximation of the pressure.

\pgfplotstableread{%
dt  Vel_l2tS_errcM  Vel_h1tS_err_M  Prs_l2tS_prsc
6.0000e-01  1.3092e-01  1.4093e-01  4.0828e-01
3.0000e-01  7.4779e-02  8.0527e-02  2.3635e-01
1.5000e-01  4.0110e-02  4.3180e-02  1.2614e-01
7.5000e-02  2.0878e-02  2.2472e-02  6.4937e-02
3.7500e-02  1.0674e-02  1.1492e-02  3.2911e-02
1.8750e-02  5.3991e-03  5.8207e-03  1.6558e-02
}\MONOexp
\pgfplotstableread{%
dt  Vel_l2tS_errcM  Vel_h1tS_err_M  Prs_l2tS_prsc
6.0000e-01  1.3091e-01  1.4092e-01  2.0381e-01
3.0000e-01  7.4775e-02  8.0522e-02  1.1687e-01
1.5000e-01  4.0109e-02  4.3177e-02  6.2258e-02
7.5000e-02  2.0878e-02  2.2470e-02  3.2186e-02
3.7500e-02  1.0674e-02  1.1490e-02  1.6386e-02
1.8750e-02  5.3990e-03  5.8192e-03  8.2711e-03
}\MONOpica
\pgfplotstableread{%
dt  Vel_l2tS_errcM  Vel_h1tS_err_M  Prs_l2tS_prsc
6.0000e-01  1.3510e-01  1.4186e-01  4.0887e-01
3.0000e-01  7.6659e-02  8.1700e-02  2.4902e-01
1.5000e-01  4.0997e-02  4.4078e-02  1.3891e-01
7.5000e-02  2.1327e-02  2.3015e-02  7.3894e-02
3.7500e-02  1.0902e-02  1.1785e-02  3.8358e-02
1.8750e-02  5.5140e-03  5.9719e-03  1.9737e-02
}\ACexp
\pgfplotstableread{%
dt  Vel_l2tS_errcM  Vel_h1tS_err_M  Prs_l2tS_prsc
6.0000e-01  1.3091e-01  1.4092e-01  2.0380e-01
3.0000e-01  7.4775e-02  8.0522e-02  1.1687e-01
1.5000e-01  4.0109e-02  4.3177e-02  6.2260e-02
7.5000e-02  2.0878e-02  2.2471e-02  3.2196e-02
3.7500e-02  1.0674e-02  1.1491e-02  1.6404e-02
1.8750e-02  5.3999e-03  5.8204e-03  8.2849e-03
}\ACpica
\pgfplotstableread{%
dt  Vel_l2tS_errcM  Vel_h1tS_err_M  Prs_l2tS_prsc
6.0000e-01  1.1310e-01  1.2232e-01  1.8105e-01
3.0000e-01  4.9324e-02  5.3612e-02  8.1222e-02
1.5000e-01  1.7689e-02  1.9334e-02  2.9896e-02
7.5000e-02  5.5112e-03  6.0800e-03  9.6885e-03
3.7500e-02  1.5629e-03  1.7690e-03  2.8884e-03
1.8750e-02  4.1774e-04  5.7504e-04  8.1212e-04
}\BDFpica
\pgfplotstableread{%
dt  Vel_l2tS_errcM  Vel_h1tS_err_M  Prs_l2tS_prsc
6.0000e-01  1.1311e-01  1.2234e-01  3.9305e-01
3.0000e-01  4.9328e-02  5.3615e-02  2.0634e-01
1.5000e-01  1.7690e-02  1.9337e-02  9.1773e-02
7.5000e-02  5.5115e-03  6.0824e-03  3.6543e-02
3.7500e-02  1.5630e-03  1.7712e-03  1.3676e-02
1.8750e-02  4.1771e-04  5.7679e-04  4.9436e-03
}\BDFexp
\pgfplotstableread{%
dt  Vel_l2tS_errcM  Vel_h1tS_err_M  Prs_l2tS_prsc
6.0000e-01  1.1804e-01  1.2365e-01  4.2749e-01
3.0000e-01  5.1715e-02  5.5231e-02  2.3195e-01
1.5000e-01  1.8786e-02  2.0615e-02  1.0703e-01
7.5000e-02  6.0211e-03  6.8198e-03  4.6700e-02
3.7500e-02  1.8393e-03  2.1973e-03  2.1411e-02
1.8750e-02  6.0036e-04  8.3204e-04  1.1489e-02
}\BOOTexp

Additional computations with higher Reynolds numbers are performed to study 
the dependence of the results delivered by the AC method on the parameter
$\ACprm$. The results are reported in \cref{tab:tgv_eta_val}. As expected, the
predictions of the AC method get closer to those of the monolithic approach
when larger values of $\ACprm$ are used. However, using larger values also
leads to higher computational costs because it affects the conditioning of the
linear system through the presence of the grad-div term. The results altogether
indicate that an appropriate choice is to set  $\ACprm=10\Rey$. Indeed, this
choice leads to a satisfactory trade-off between the accuracy and the
conditioning of the linear systems.

\begin{table}
  \centering
  \caption{2D Taylor--Green vortex. Space-time velocity and pressure errors.
    Explicit convection. AC method with $\ACprm\in\Set{\Rey,10\Rey,100\Rey}$.
    The errors obtained with the monolithic approach are given in parenthesis
    as reference. First group of columns: $\Rey\approx33$, Cartesian mesh
    composed of $128^2$ cells, $T=40$, $\deltat=\frac{T}{64}=0.625$. Second
    group of columns: $\Rey=100$, Cartesian mesh composed of $512^2$ cells,
  $T=120$, $\deltat=\frac{T}{64}=1.875$.}
  \label{tab:tgv_eta_val}
  \begingroup
  \scriptsize
  \begin{tabular}{c|ccc:c|ccc:c}
    & \multicolumn{4}{c|}{$\Rey\approx33$} & \multicolumn{4}{c}{$\Rey=100$}\\
    $\ACprm$ & $\Rey$ & $10\Rey$ & $100\Rey$ & (MONO) & $\Rey$ & $10\Rey$ & $100\Rey$ & (MONO)\\
    \hline
    $\ell^2(\Lebesgue2)$ velocity 
    & \pgfmathprintnumber[prnt]{1.0530e-02}
    & \pgfmathprintnumber[prnt]{1.4097e-03}
    & \pgfmathprintnumber[prnt]{1.2436e-03}
    &(\pgfmathprintnumber[prnt]{1.2369e-03})
    & \pgfmathprintnumber[prnt]{2.9127e-02}
    & \pgfmathprintnumber[prnt]{3.9177e-03}
    & \pgfmathprintnumber[prnt]{1.3821e-03}
    &(\pgfmathprintnumber[prnt]{1.1357e-03})\\
    $\ell^2(\Hilbert1)$ velocity 
    & \pgfmathprintnumber[prnt]{1.2584e-02}
    & \pgfmathprintnumber[prnt]{1.0391e-02}
    & \pgfmathprintnumber[prnt]{1.0328e-02}
    &(\pgfmathprintnumber[prnt]{1.0320e-02})
    & \pgfmathprintnumber[prnt]{2.0143e-02}
    & \pgfmathprintnumber[prnt]{8.2460e-03}
    & \pgfmathprintnumber[prnt]{7.8308e-03}
    &(\pgfmathprintnumber[prnt]{7.8024e-03})\\
    $\ell^2(\Lebesgue2)$ pressure 
    & \pgfmathprintnumber[prnt]{2.8050e-02}
    & \pgfmathprintnumber[prnt]{1.1862e-02}
    & \pgfmathprintnumber[prnt]{1.0313e-02}
    &(\pgfmathprintnumber[prnt]{9.9976e-03})
    & \pgfmathprintnumber[prnt]{6.0206e-02}
    & \pgfmathprintnumber[prnt]{1.6859e-02}
    & \pgfmathprintnumber[prnt]{1.2464e-02}
    &(\pgfmathprintnumber[prnt]{1.1885e-02})
  \end{tabular}
  \endgroup
\end{table}

\subsubsection{CFL conditions with explicit convection}
\label{ssub:tgv_stab}
We investigate numerically the stability of the monolithic and AC approaches when an
explicit treatment of the convection term is used. For these tests, we consider
\begin{enumerate*}[label=(\roman*)]
  \item three Reynolds numbers, $\Rey\in\Set{200,500,700}$,
  \item an observation time $T$ such that $T\,\Rey=10^4$,
  \item $\ACprm=10\,\Rey$ whenever the AC method is considered.
\end{enumerate*} 
We flag a computation as having diverged if, for some $n\geq1$, we have
\begin{equation}
  \label{eq:TGV_div_crit}
  \EnrKin(\hVex\uh^n)>1.1\,\EnrKin(\hVex\uh^0)=1.1\,\EnrKin(\hvprj\uh(\Vex_0))\,,
\end{equation}
where the definition of the discrete kinetic energy $\EnrKin(\cdot)$ is given in
\eqref{eq:kin_en}. Notice that the solution \eqref{eq:TGV} goes exponentially
towards $0$ with respect to time. Thus, a failure to satisfy \eqref{eq:TGV_div_crit}
is a symptom of stability issues. 
We denote by $\dtc$ the critical time-step, that is,
the largest $\deltat$ ensuring that the computation does not diverge. 
We seek a resolution of $1\%$, meaning that the gap between $\dtc$ and the
smallest $\deltat$ leading to divergence is less than $1\%$ of $\dtc$.
Whenever $\deltat>\dtc$, we define the divergence time, $\Tdiv$, as the smallest $t^n$
for which \eqref{eq:TGV_div_crit} is satisfied.

\begin{figure}
  \centering
  \includegraphics[width=0.4\textwidth]{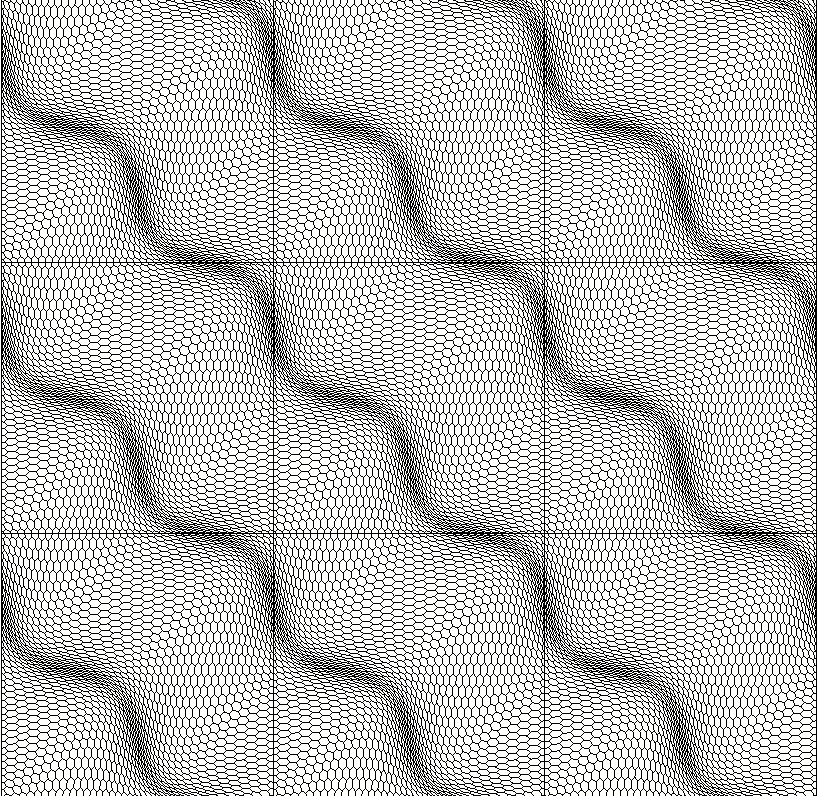}
  \caption{2D polygonal mesh used for CFL study.}
  \label{fig:2D_prg}
\end{figure}

\begin{figure}
  \centering
  \begin{tikzpicture}
    \begin{groupplot}[
      group style={group size=2 by 1,group name=G,horizontal sep=0.5cm,
        ylabels at=edge left,y descriptions at=edge left,},
      width=0.45\linewidth,base plot,xlabel={$\deltat\,\Rey$},
      xmin=5,xmax=6.5,ymin=1200,ymax=4000,
      xtick distance=0.5,ytick distance=1000,
      ]
      \nextgroupplot[title={Cartesian mesh \ref{mrk:stbMNfst}},ylabel={$\Tdiv\,\Rey$}]
        \addplot[gener,solid]  table[x expr={\thisrow{R2_dt}*200},  y expr={\thisrow{R2_T}*200}]  {\failedtable};\label{mrk:stbtwoh}
        \addplot[gener,dashed] table[x expr={\thisrow{R5_dt}*500},  y expr={\thisrow{R5_T}*500}]  {\failedtable};\label{mrk:stbfiveh}
        \addplot[gener,dotted] table[x expr={\thisrow{R7_dt}*700},  y expr={\thisrow{R7_T}*700}]  {\failedtable};\label{mrk:stbsevenh}
        \addplot[stbMN,solid]  table[x expr={\thisrow{R2_dt}*200},  y expr={\thisrow{R2_T}*200}]  {\failedtable};\label{mrk:stbMNfst}
        \addplot[stbMN,dashed] table[x expr={\thisrow{R5_dt}*500},  y expr={\thisrow{R5_T}*500}]  {\failedtable};
        \addplot[stbMN,dotted] table[x expr={\thisrow{R7_dt}*700},  y expr={\thisrow{R7_T}*700}]  {\failedtable};
        \nextgroupplot[title={Polygonal mesh \ref{mrk:stbACfst}}]
        \addplot[stbAC,solid]  table[x expr={\thisrow{R2_dt}*200},  y expr={\thisrow{R2_T}*200}]  {\failedACtable};\label{mrk:stbACfst}
        \addplot[stbAC,dashed] table[x expr={\thisrow{R5_dt}*500},  y expr={\thisrow{R5_T}*500}]  {\failedACtable};
        \addplot[stbAC,dotted] table[x expr={\thisrow{R7_dt}*700},  y expr={\thisrow{R7_T}*700}]  {\failedACtable};
    \end{groupplot}
  \end{tikzpicture}

  \begin{tabular}{ccc}
    $\Rey=200$ \ref{mrk:stbtwoh} & $\Rey=500$ \ref{mrk:stbfiveh} & $\Rey=700$ \ref{mrk:stbsevenh} 
  \end{tabular}
  \caption{2D Taylor--Green vortex. Monolithic approach and explicit
    convection. Divergence time $\Tdiv$ (scaled by $\Rey$), for different
    choices of $\deltat$ (scaled by $\Rey$). Left: Cartesian mesh composed of 16{,}384
cells. Right: polygonal mesh composed of 15{,}129 cells.}
    \label{fig:TGV_Tdiv}
\end{figure}
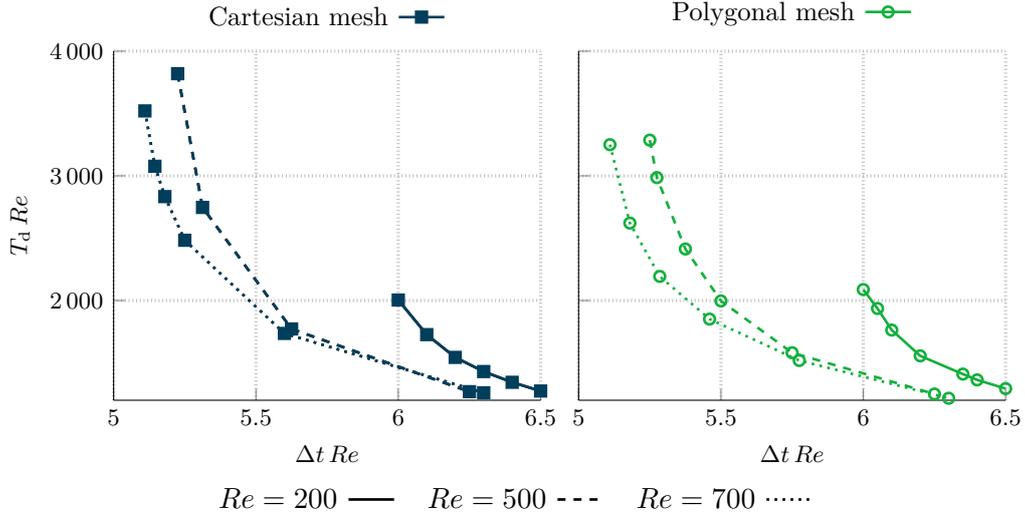

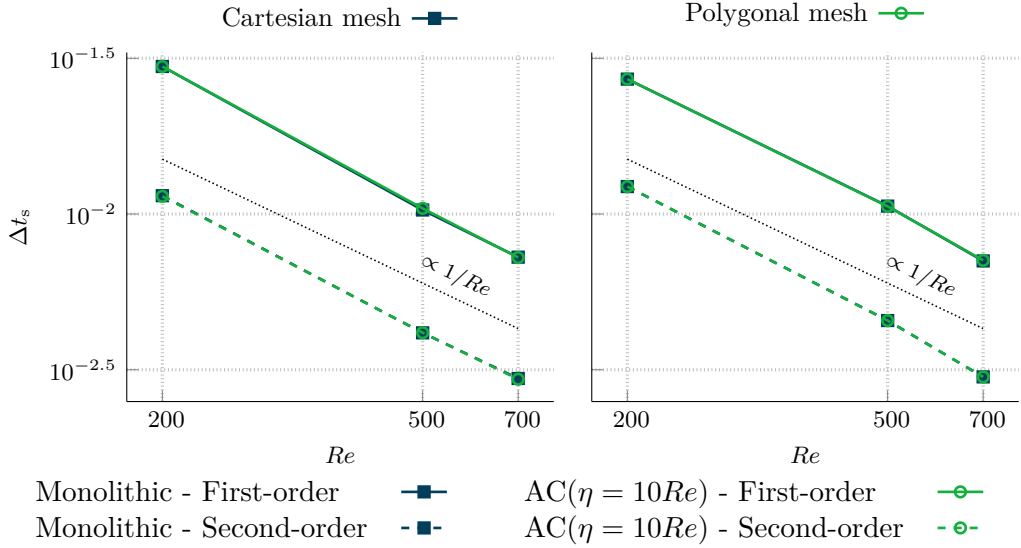
\begin{figure}
  \centering
  \begin{tikzpicture}
    \begin{groupplot}[
      group style={group size=2 by 1,group name=G,horizontal sep=0.5cm,
        ylabels at=edge left,y descriptions at=edge left,},
      width=0.45\linewidth,base plot,xmode=log,ymode=log,
      ylabel={$\dtc$},xlabel={\Rey},xtick={data},
      ymin=2.5e-3,ymax=3.3e-2,log x ticks with fixed point,
      table/x=Re,
      ]
      \nextgroupplot[title={Cartesian mesh \ref{mrk:stbMNfst}}]
        \addplot[stbMN,stbfst] table[y=1st_mn]    {\CFLtable};
        \addplot[stbAC,stbfst] table[y=1st_ac_10] {\CFLtable};
        \addplot[stbMN,stbscn] table[y=2nd_mn]    {\CFLtable};\label{mrk:stbMNscn}
        \addplot[stbAC,stbscn] table[y=2nd_ac_10] {\CFLtable};\label{mrk:stbACscn}
        \addplot[base,order]  expression[domain=200:700] {3/x} node[ordernode] {$\propto 1/\Rey$};
      \nextgroupplot[title={Polygonal mesh \ref{mrk:stbACfst}}]
        \addplot[stbMN,stbfst] table[y=1st_mn]    {\CFLtablePoly};
        \addplot[stbAC,stbfst] table[y=1st_ac_10] {\CFLtablePoly};
        \addplot[stbMN,stbscn] table[y=2nd_mn]    {\CFLtablePoly};
        \addplot[stbAC,stbscn] table[y=2nd_ac_10] {\CFLtablePoly};
        \addplot[base,order]  expression[domain=200:700] {3/x} node[ordernode] {$\propto 1/\Rey$};
    \end{groupplot}
  \end{tikzpicture}
  \begin{tabular}{ll@{\hspace{1cm}}ll}
    Monolithic - First-order & \ref{mrk:stbMNfst} & AC($\ACprm=10\Rey$) - First-order & \ref{mrk:stbACfst} \\
    Monolithic - Second-order & \ref{mrk:stbMNscn} & AC($\ACprm=10\Rey$) - Second-order & \ref{mrk:stbACscn}
  \end{tabular}
  \caption{2D Taylor--Green vortex. Critical time-step for stability, $\dtc$
    (up to $1\%$ resolution), as a function of the Reynolds number, $\Rey$, for several 
  time-stepping schemes. Left: Cartesian mesh composed of 16{,}384
cells. Right: polygonal mesh composed of 15{,}129 cells.}
  \label{fig:TGV_stb}
\end{figure}

Two meshes have been considered, a Cartesian mesh and a polygonal mesh (illustrated in
\cref{fig:2D_prg}), which have, respectively, 16{,}384 and 15{,}129 cells. In
\cref{fig:TGV_Tdiv}, we plot the divergence time, $\Tdiv$, for
different time-step values. We observe that the divergence time sharply
increases as the time step is decreased towards the critical value. This confirms
that using an observation time such that $T\,\Rey=10^4$ is reasonable to flag
stable/unstable computations. \Cref{fig:TGV_stb} shows the critical time-step values
for the various time-stepping schemes and the two meshes. For each case, a
dependence of the critical
time-step on the inverse of the Reynolds number is observed. 
For convenience, the results are also reported in
\cref{tab:TGV_sum_fst_vs_scn,tab:TGV_sum_fst_vs_scn_poly} for 
the Cartesian and
the polygonal meshes, respectively, including a quantitative 
comparison between the first- and the 
second-order schemes. Generally speaking, the AC method turns out to be as stable as the
monolithic scheme, and no significant differences on the critical time-step are
observed. Moreover, the first-order schemes are more stable than the second-order ones,
and allow one to
choose time-step values more than two times larger than those possible with
second-order schemes. Finally, the use of a polygonal mesh only marginally 
affects the value of the critical time-step.

\pgfplotstableread{%
R2_dt   R2_T   R5_dt     R5_T   R7_dt     R7_T
3.00e-2 10.02  1.0450e-2 7.6389 7.3000e-3 5.0297
3.05e-2 8.6315 1.0625e-2 5.4931 7.3500e-3 4.3953
3.10e-2 7.719  1.1250e-2 3.5437 7.4000e-3 4.0478
3.15e-2 7.1505 1.2500e-2 2.5375 7.5000e-3 3.5475
3.20e-2 6.72   NaN       NaN    8.0000e-3 2.4800
3.25e-2 6.37   NaN       NaN    9.0000e-3 1.8000
}\failedtable
\pgfplotstableread{%
R2_dt       R2_T     R5_dt       R5_T     R7_dt       R7_T
3.0000e-2   10.440   1.0500e-2   6.5730   7.3000e-3   4.6428
3.0250e-2   9.6800   1.0550e-2   5.9713   7.4000e-3   3.7444
3.0500e-2   8.8145   1.0750e-2   4.8267   7.5500e-3   3.1333
3.1000e-2   7.7810   1.1000e-2   3.9930   7.8000e-3   2.6442
3.1750e-2   7.0485   1.1500e-2   3.1625   8.2500e-3   2.1698
3.2000e-2   6.8160   1.2500e-2   2.5000   9.0000e-3   1.7370
3.2500e-2   6.4675   NaN         NaN      NaN         NaN   
}\failedACtable
\pgfplotstableread{%
Re   1st_mn    2nd_mn   1st_ac_10 2nd_ac_10
200  2.9750e-2 1.145e-2 2.9750e-2 1.1400e-2
500  1.03e-2   4.155e-3 1.0425e-2 4.1600e-3
700  7.27e-3   2.965e-3 7.2500e-3 2.9500e-3
}\CFLtable
\pgfplotstableread{%
Re   1st_mn    2nd_mn   1st_ac_10  2nd_ac_10
200  2.71e-2   1.225e-2 2.71e-2    1.225e-2 
500  1.06e-2   4.55e-3  1.058e-2   4.55e-3
700  7.07e-3   3.00e-3  7.10e-3    3.00e-3
}\CFLtablePoly

\begin{table}
  \centering
  \caption{2D Taylor--Green vortex. Critical time-step for
    stability, $\dtc$, as a function of the Reynolds number, $\Rey$. Comparison of
    first- and second-order schemes. Cartesian mesh composed of 16{,}384
cells.}
  \label{tab:TGV_sum_fst_vs_scn}
  \begin {tabular}{c|ccc|ccc}%
  & \multicolumn {3}{c|}{Monolithic} & \multicolumn {3}{c}{AC($\ACprm =10\Rey $)}\\$\Rey $&$1^{\mathrm {st}}$&$2^{\mathrm {nd}}$&$\frac {1^{\rm {}st}}{2^{\rm {}nd}}$&$1^{\mathrm {st}}$&$2^{\mathrm {nd}}$&$\frac {1^{\rm {}st}}{2^{\rm {}nd}}$\\\hline %
  \pgfutilensuremath {200}&\pgfutilensuremath {2.98e{-}2}&\pgfutilensuremath {1.15e{-}2}&\pgfutilensuremath {2.60}&\pgfutilensuremath {2.98e{-}2}&\pgfutilensuremath {1.14e{-}2}&\pgfutilensuremath {2.61}\\%
  \pgfutilensuremath {500}&\pgfutilensuremath {1.03e{-}2}&\pgfutilensuremath {4.16e{-}3}&\pgfutilensuremath {2.48}&\pgfutilensuremath {1.04e{-}2}&\pgfutilensuremath {4.16e{-}3}&\pgfutilensuremath {2.51}\\%
  \pgfutilensuremath {700}&\pgfutilensuremath {7.27e{-}3}&\pgfutilensuremath {2.97e{-}3}&\pgfutilensuremath {2.45}&\pgfutilensuremath {7.25e{-}3}&\pgfutilensuremath {2.95e{-}3}&\pgfutilensuremath {2.46}\\%
  \end {tabular}%
\end{table}
\begin{table}
  \centering
  \caption{2D Taylor--Green vortex. Critical time-step for
    stability, $\dtc$, as a function of the Reynolds number, $\Rey$. Comparison of
    first- and second-order schemes. Polygonal mesh composed of 15{,}129 cells.}
  \label{tab:TGV_sum_fst_vs_scn_poly}
  \begin {tabular}{c|ccc|ccc}%
  & \multicolumn {3}{c|}{Monolithic} & \multicolumn {3}{c}{AC($\ACprm =10\Rey $)}\\$\Rey $&$1^{\mathrm {st}}$&$2^{\mathrm {nd}}$&$\frac {1^{\rm {}st}}{2^{\rm {}nd}}$&$1^{\mathrm {st}}$&$2^{\mathrm {nd}}$&$\frac {1^{\rm {}st}}{2^{\rm {}nd}}$\\\hline %
  \pgfutilensuremath {200}&\pgfutilensuremath {2.71e{-}2}&\pgfutilensuremath {1.23e{-}2}&\pgfutilensuremath {2.21}&\pgfutilensuremath {2.71e{-}2}&\pgfutilensuremath {1.23e{-}2}&\pgfutilensuremath {2.21}\\%
  \pgfutilensuremath {500}&\pgfutilensuremath {1.06e{-}2}&\pgfutilensuremath {4.55e{-}3}&\pgfutilensuremath {2.33}&\pgfutilensuremath {1.06e{-}2}&\pgfutilensuremath {4.55e{-}3}&\pgfutilensuremath {2.33}\\%
  \pgfutilensuremath {700}&\pgfutilensuremath {7.07e{-}3}&\pgfutilensuremath {3.00e{-}3}&\pgfutilensuremath {2.36}&\pgfutilensuremath {7.10e{-}3}&\pgfutilensuremath {3.00e{-}3}&\pgfutilensuremath {2.37}\\%
  \end {tabular}%
\end{table}

\subsection{3D modified Taylor--Green vortex}
\label{sec:3D}

\begin{figure}
  \centering
  \includegraphics[width=.3\linewidth]{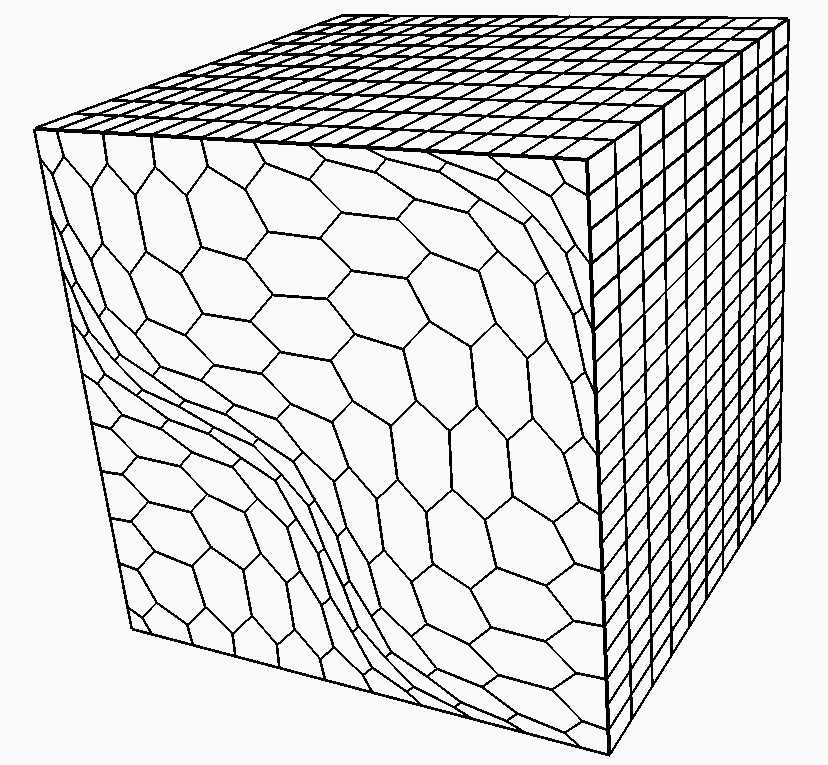}
  \caption{Coarsest mesh of the PrG sequence of \cite{FVCA6} composed of prisms with
  polygonal basis.}
  \label{fig:prg}
\end{figure}

We consider the 3D modified Taylor--Green vortex solution of \cite[Benchmark case
2.2]{FVCA8} and adapt it so that it is a solution to the Navier--Stokes equations. In
particular, we make the solution time-dependent by considering a sinusoidal amplitude. The
exact solution reads as follows:
\begin{equation}
  \left\{
  \begin{aligned}
    \Vex(x,y,z) \coloneqq{}& \alpha(t) \Vex^{\prime}(x,y,z)\,, \\
    \Pex(x,y,z) \coloneqq{}& \alpha(t) \Pex^{\prime}(x,y,z)\,, \\
    \alpha(t) \coloneqq{}& \sin(8\pi t)\,,\\
    \Vex^{\prime}(x,y,z) \coloneqq{}&\left[ %
      \begin{aligned}
        -2&\cos(2\pi x)\sin(2\pi y) \sin(2\pi z) \\&\sin(2\pi x)\cos(2\pi
        y) \sin(2\pi z) \\&\sin(2\pi x)\sin(2\pi y) \cos(2\pi z)
      \end{aligned} \right]\,, \\
    \Pex^{\prime}(x,y,z) \coloneqq{}& -6\pi \sin(2\pi x)\sin(2\pi y) \sin(2\pi z)\,,
  \end{aligned}
  \right.
\end{equation}
and the resulting source term is such that
\begin{equation}
  \left\{
  \begin{aligned}
    \Foex(x,y,z) \coloneqq{}& \alpha(t)\Foex^{\prime}(x,y,z) + 8\,\pi\cos(8\,\pi t)\Vex^{\prime}(x,y,z)+\\
                            &-\frac{\alpha^2(t)}{4} \left[ %
      \begin{aligned}
        -2&\sin(4\pi x)(\cos(4\pi y)+\cos(4\pi z)-2)\\
          &\sin(4\pi y)(\cos(4\pi x)-2\cos(4\pi z)+1)\\
          &\sin(4\pi z)(\cos(4\pi x)-2\cos(4\pi y)+1)
      \end{aligned} \right]\,,\\
    \Foex^{\prime}(x,y,z) \coloneqq{}& [-36\pi^2\cos(2\pi x)\sin(2\pi y) \sin(2\pi z),\,0,\,0]^T\,.
  \end{aligned}
  \right.
\end{equation}
The domain is the unit cube $\Dom\coloneqq (0,1)^3$ and the viscosity is $\Visc\coloneqq1$.
Hence, since it is again reasonable to set 
$L\coloneqq1$ and $U\coloneqq1$ for the reference length and velocity, respectively, 
we obtain a Reynolds number of
$\Rey=1$. The observation time is set to $T\coloneqq 2$ and the considered 
time-step values are
$\deltat\in\Set{\frac{T}{32},\,\frac{T}{64},\,\frac{T}{128},\,\frac{T}{256}}$. The
convection term is always treated explicitly.
A quite refined
mesh composed of prisms with polygonal bases is considered.
Namely, the mesh obtained from the PrG sequence proposed
in the benchmark \cite{FVCA6} by gluing six of the finest meshes side by side 
and then rescaling the resulting mesh to the desired dimension
(an example is shown in \cref{fig:prg}). The resulting mesh
has more than 14M cells and leads to a final coupled system of more than 170M
DoFs (after static condensation). 

Only iterative solvers are considered in this 3D test case.
In particular, the linear systems obtained with the AC method
are solved by means of a Jacobi-preconditioned Conjugate Gradient (CG) solver. In
order to deal with the saddle point systems obtained with the monolithic
approach, an iterative Golub-Kahan Bidiagonalization (GKB) \cite{Ario13,AKRT18}
procedure without augmentation is chosen; the internal solver is a Conjugate
Gradient with a K-cycle Algebraic Multigrdid preconditioner~\cite{NoVa08}.
All these solvers are natively available in \CS{}. 
The computations have been run on the GAIA
cluster of EDF\footnote{243\textsuperscript{rd} of the TOP500 list of November
2020} on up to 700 cores. Two values of the AC parameter have been considered:
$\ACprm=10\Rey$ and $\ACprm=50\Rey$. 

\cref{fig:MTGV_cvg} reports the convergence results for the normalized space-time
velocity and pressure errors as defined in \eqref{eq:def_errors_num}.
These results corroborate the optimal convergence in time for both
first- and second-order schemes. A slight decrease of the convergence
rate for the velocity $\ell^2(\Hilbert1)$-error is observed for the 
finest $\deltat$ and the second-order schemes, which indicates that the space discretization error tends to become dominant.  
Moreover, we remark that the higher the value of $\ACprm$, the more accurate the results
obtained with the AC method.
Actually, for $\ACprm=50\Rey$, the errors are essentially 
superimposed to those obtained with the
monolithic approach. 

\pgfmathsetmacro{\LLVel}{sqrt(0.75)}
\pgfmathsetmacro{\LHVel}{3*(4*atan(1))}
\pgfmathsetmacro{\LLPrs}{3*(4*atan(1))*.5}

\pgfplotstableread{%
dt  Vel_l2tS_errcM  Vel_h1tS_err_M  Prs_l2tS_prsc  Tot_time
6.2500e-02  8.5549e-02  5.1294e-01  1.1352e+00  4.7782e+06
3.1250e-02  4.7150e-02  2.8274e-01  5.9095e-01  9.9094e+06
1.5625e-02  2.4781e-02  1.5100e-01  3.0565e-01  1.9963e+07
7.8125e-03  1.2724e-02  8.2541e-02  1.5502e-01  4.2289e+07
}\monofst
\pgfplotstableread{%
dt  Vel_l2tS_errcM  Vel_h1tS_err_M  Prs_l2tS_prsc  Tot_time
6.2500e-02  1.3060e-01  7.8726e-01  2.1624e+00  1.2306e+06
3.1250e-02  8.1976e-02  4.8489e-01  1.4569e+00  2.2456e+06
1.5625e-02  4.5041e-02  2.6526e-01  8.0893e-01  3.9105e+06
7.8125e-03  2.3117e-02  1.3870e-01  4.1071e-01  6.4318e+06
}\actfst
\pgfplotstableread{%
dt  Vel_l2tS_errcM  Vel_h1tS_err_M  Prs_l2tS_prsc  Tot_time
6.2500e-02  8.9285e-02  5.3528e-01  1.2111e+00  2.4790e+06
3.1250e-02  4.9414e-02  2.9561e-01  6.5754e-01  4.4846e+06
1.5625e-02  2.5971e-02  1.5761e-01  3.4167e-01  7.7284e+06
7.8125e-03  1.3327e-02  8.5672e-02  1.7260e-01  1.2865e+07
}\achfst
\pgfplotstableread{%
dt  Vel_l2tS_errcM  Vel_h1tS_err_M  Prs_l2tS_prsc  Tot_time
6.2500e-02  8.0169e-02  5.0949e-01  1.8168e+00  7.2807e+06
3.1250e-02  2.5418e-02  1.6273e-01  5.2148e-01  1.4152e+07
1.5625e-02  6.7896e-03  5.3773e-02  1.3209e-01  3.0120e+07
7.8125e-03  1.7349e-03  3.4872e-02  3.3245e-02  6.3394e+07
}\monoscn
\pgfplotstableread{%
dt  Vel_l2tS_errcM  Vel_h1tS_err_M  Prs_l2tS_prsc  Tot_time
6.2500e-02  9.0774e-02  5.7376e-01  1.9684e+00  2.6825e+06
3.1250e-02  3.3315e-02  2.0789e-01  7.4169e-01  4.7827e+06
1.5625e-02  9.8811e-03  6.8991e-02  2.2987e-01  8.3948e+06
7.8125e-03  2.6950e-03  3.7095e-02  6.9687e-02  1.3727e+07
}\actscn
\pgfplotstableread{%
dt  Vel_l2tS_errcM  Vel_h1tS_err_M  Prs_l2tS_prsc  Tot_time
6.2500e-02  8.0615e-02  5.1177e-01  1.8064e+00  5.5209e+06
3.1250e-02  2.5623e-02  1.6369e-01  5.3531e-01  1.0062e+07
1.5625e-02  6.8709e-03  5.4095e-02  1.3952e-01  1.7147e+07
7.8125e-03  1.7652e-03  3.4925e-02  3.7572e-02  2.8435e+07
}\achscn
\newcommand{\domain}{7.125e-3:6.25e-2}
\begin{figure}
  \centering
  \begin{tikzpicture}
    \begin{groupplot}[
      group style={group size=3 by 1,},
      xmode=log,ymode=log,base plot,width=0.365\linewidth,
      xticklabels={$\frac{T}{256}$,$\frac{T}{128}$,$\frac{T}{64}$,$\frac{T}{32}$,$\frac{T}{16}$,$\frac{T}{8}$,$\frac{T}{4}$},
      xtick={7.8125e-3,1.5625e-2,3.125e-2,6.25e-2,0.125,0.25,0.5},
      ]
    \nextgroupplot[
        title={Velocity - $\ell^2(\Lebesgue2)$},
        table/x=dt,table/y expr={\thisrow{Vel_l2tS_errcM}/\LLVel},
      ]
      \addplot[fst,order]  expression[domain=\domain] {6*x}    node[ordernode] {1\fst{} ord};
      \addplot[scn,order]  expression[domain=\domain] {20*x*x} node[ordernode,pos=0.2,below] {2\scn{} ord};
      \addplot[fst,gener]      table {\monofst};\label{mrk:fst}  
      \addplot[fst,mono,pzero] table {\monofst};\label{mrk:mono}
      \addplot[fst,ac,pten]    table {\actfst}; \label{mrk:act}
      \addplot[fst,ac,phun]    table {\achfst}; \label{mrk:ach}
      \addplot[scn,gener]      table {\monoscn};\label{mrk:scn}  
      \addplot[scn,mono,pzero] table {\monoscn};
      \addplot[scn,ac,pten]    table {\actscn};
      \addplot[scn,ac,phun]    table {\achscn};
    \nextgroupplot[
        title={Velocity - $\ell^2(\Hilbert1)$},
        table/x=dt,table/y expr={\thisrow{Vel_h1tS_err_M}/\LHVel},
      ]
      \addplot[fst,order]  expression[domain=\domain] {.06*x}    node[ordernode] {1\fst{} ord};
      \addplot[scn,order]  expression[domain=\domain] {.2*x*x} node[ordernode,pos=0.2,below] {2\scn{} ord};
      \addplot[fst,mono,pzero] table {\monofst};
      \addplot[fst,ac,pten]    table {\actfst};
      \addplot[fst,ac,phun]    table {\achfst};
      \addplot[scn,mono,pzero] table {\monoscn};
      \addplot[scn,ac,pten]    table {\actscn};
      \addplot[scn,ac,phun]    table {\achscn};
    \nextgroupplot[
        title={Pressure - $\ell^2(\Lebesgue2)$},
        table/x=dt,table/y expr={\thisrow{Prs_l2tS_prsc}/\LLPrs},
      ]
      \addplot[fst,order]  expression[domain=\domain] {.25*x}    node[ordernode] {1\fst{} ord};
      \addplot[scn,order]  expression[domain=\domain] {1*x*x} node[ordernode,pos=0.2,below] {2\scn{} ord};
      \addplot[fst,mono,pzero] table {\monofst};
      \addplot[fst,ac,pten]    table {\actfst};
      \addplot[fst,ac,phun]    table {\achfst};
      \addplot[scn,mono,pzero] table {\monoscn};
      \addplot[scn,ac,pten]    table {\actscn};
      \addplot[scn,ac,phun]    table {\achscn};
    \end{groupplot}
  \end{tikzpicture}
  \begin{tabular}{r@{: }lll}
    Solver     &  Monolithic \ref{mrk:mono} & AC($\ACprm=10$) \ref{mrk:act} & AC($\ACprm=50$) \ref{mrk:ach} \\
    Time order & 1\fst{} \ref{mrk:fst} & 2\scn{} \ref{mrk:scn}
  \end{tabular}
  \caption{3D modified Taylor--Green vortex. Convergence results for the space-time velocity and pressure errors. $\Rey=1$.}
  \label{fig:MTGV_cvg}
\end{figure}
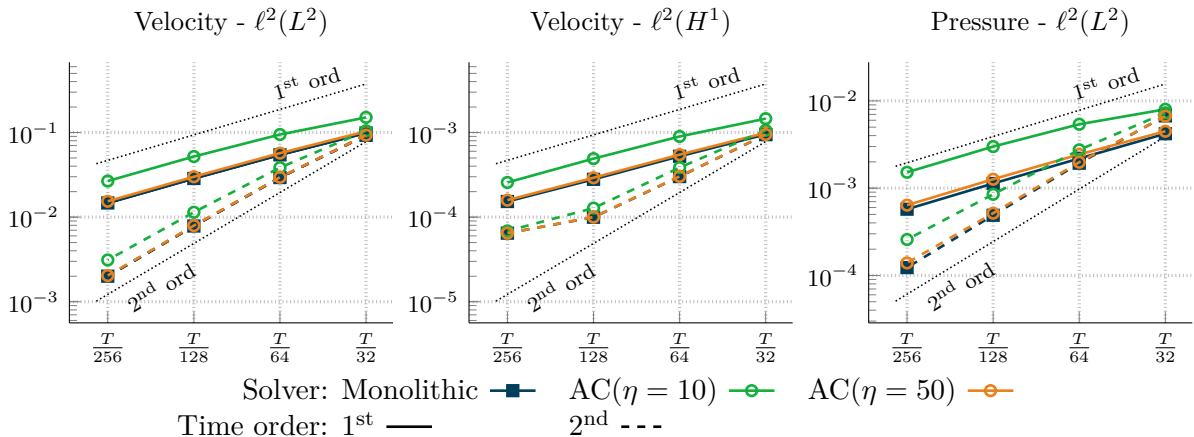

Our aim is now to compare the computational efficiency of the monolithic
and AC methods by studying the achieved velocity and pressure errors
as a function of the computational cost (defined as the product
between the elapsed time and the number of computing cores). 
As a preliminary step, 
computations are run in order to find reasonably optimized values for the tolerances of
the iterative solvers. The results for the first-order AC and monolithic
approaches are, respectively, gathered in
\cref{tab:mtgv_tol_ac,tab:mtgv_tol_mono}. A tolerance of $10^{-4}$ appears to be
a reasonable choice for the CG iterative solver in the AC context, for both considered
values of $\ACprm$. In order to anticipate the greater accuracy of second-order schemes,
the tolerance is decreased to $10^{-5}$ for the second-order AC time-stepping scheme. 
Similarly, for the monolithic approach, the tolerances for the GKB and CG iterative solvers are both set to $10^{-4}$ and $10^{-5}$, respectively, for the first- and second-order time-stepping schemes.

\begin{table}
  \centering
  \caption{3D modified Taylor--Green vortex. Space-time errors for different
    tolerances of the CG linear solver. First-order AC time-stepping scheme, $\deltat=\frac{T}{128}$. 
    The adopted tolerances for further calculations are circled.}
  \label{tab:mtgv_tol_ac}
  \begin {tabular}{c|ccc}%
  $\epsilon $(CG)&$\ell ^2(\Lebesgue 2)~\Vex $&$\ell ^2(\Hilbert 1)~\Vex $&$\ell ^2(\Lebesgue 2)~\Pex $\\\hline &\multicolumn {3}{c}{$\ACprm =10$}\\%
  \pgfutilensuremath {1e{-}2}&\pgfutilensuremath {6.58e{-}2}&\pgfutilensuremath {6.76e{-}4}&\pgfutilensuremath {3.17e{-}3}\\%
  \tikz[baseline=(A.base)]\node[draw=black,rectangle,rounded corners,](A){\pgfutilensuremath {1e{-}4}};&\pgfutilensuremath {5.20e{-}2}&\pgfutilensuremath {4.91e{-}4}&\pgfutilensuremath {3.00e{-}3}\\%
  \pgfutilensuremath {1e{-}6}&\pgfutilensuremath {5.20e{-}2}&\pgfutilensuremath {4.91e{-}4}&\pgfutilensuremath {3.00e{-}3}\\\hline &\multicolumn {3}{c}{$\ACprm =50$}\\%
  \pgfutilensuremath {1e{-}2}&\pgfutilensuremath {3.83e{-}2}&\pgfutilensuremath {4.22e{-}4}&\pgfutilensuremath {1.49e{-}3}\\%
  \tikz[baseline=(A.base)]\node[draw=black,rectangle,rounded corners,](A){\pgfutilensuremath {1e{-}4}};&\pgfutilensuremath {3.00e{-}2}&\pgfutilensuremath {2.92e{-}4}&\pgfutilensuremath {1.27e{-}3}\\%
  \pgfutilensuremath {1e{-}6}&\pgfutilensuremath {3.00e{-}2}&\pgfutilensuremath {2.92e{-}4}&\pgfutilensuremath {1.27e{-}3}\\%
  \end {tabular}%
\end{table}
\begin{table}
  \centering
  \caption{3D modified Taylor--Green vortex. Space-time errors for different
    tolerances of the GKB and CG linear solvers. First-order monolithic time-stepping scheme,
    $\deltat=\frac{T}{128}$. 
    The adopted tolerances for further calculations are circled.}
  \label{tab:mtgv_tol_mono}
  \begin {tabular}{cc|ccc}%
  $\epsilon $(GKB)&$\epsilon $(CG)&$\ell ^2(\Lebesgue 2)~\Vex $&$\ell ^2(\Hilbert 1)~\Vex $&$\ell ^2(\Lebesgue 2)~\Pex $\\\hline %
  \multirow{3}{*}{\pgfutilensuremath {1e{-}2}}
    &\pgfutilensuremath {1e{-}2}&\pgfutilensuremath {4.27e{-}2}&\pgfutilensuremath {8.31e{-}4}&\pgfutilensuremath {2.80e{-}3}\\%
    &\pgfutilensuremath {1e{-}4}&\pgfutilensuremath {2.84e{-}2}&\pgfutilensuremath {3.21e{-}4}&\pgfutilensuremath {1.28e{-}3}\\%
    &\pgfutilensuremath {1e{-}6}&\pgfutilensuremath {2.83e{-}2}&\pgfutilensuremath {3.20e{-}4}&\pgfutilensuremath {1.28e{-}3}\\\hline %
  \multirow{2}{*}{ \tikz[baseline=(A.base)]\node[draw=black,rectangle,rounded corners,](A){\pgfutilensuremath {1e{-}4}};}
    &\tikz[baseline=(A.base)]\node[draw=black,rectangle,rounded corners,](A){\pgfutilensuremath {1e{-}4}};&\pgfutilensuremath {2.86e{-}2}&\pgfutilensuremath {2.80e{-}4}&\pgfutilensuremath {1.13e{-}3}\\%
    &\pgfutilensuremath {1e{-}6}&\pgfutilensuremath {2.87e{-}2}&\pgfutilensuremath {2.80e{-}4}&\pgfutilensuremath {1.13e{-}3}\\\hline %
  \multirow{2}{*}{\pgfutilensuremath {1e{-}6}}
    &\pgfutilensuremath {1e{-}4}&\pgfutilensuremath {2.86e{-}2}&\pgfutilensuremath {2.80e{-}4}&\pgfutilensuremath {1.13e{-}3}\\%
    &\pgfutilensuremath {1e{-}6}&\pgfutilensuremath {2.87e{-}2}&\pgfutilensuremath {2.80e{-}4}&\pgfutilensuremath {1.13e{-}3}\\%
  \end {tabular}%
\end{table}

\begin{figure}
  \centering
  \begin{tikzpicture}
    \begin{groupplot}[
      group style={group size=3 by 1,group name=G,},
      xmode=log,ymode=log,base plot,width=0.365\linewidth,
      ]
    \nextgroupplot[
        xlabel={Elapsed time $\times$ cores [\si{s}]},
        title={Velocity - $\ell^2(\Lebesgue2)$},
        table/x=Tot_time,table/y expr={\thisrow{Vel_l2tS_errcM}/\LLVel},
      ]
      \addplot[fst,mono,pzero] table {\monofst};
      \addplot[fst,ac,pten]    table {\actfst};
      \addplot[fst,ac,phun]    table {\achfst};
      \addplot[scn,mono,pzero] table {\monoscn};
      \addplot[scn,ac,pten]    table {\actscn};
      \addplot[scn,ac,phun]    table {\achscn};
    \nextgroupplot[
        xlabel={Elapsed time $\times$ cores [\si{s}]},
        title={Velocity - $\ell^2(\Hilbert1)$},
        table/x=Tot_time,table/y expr={\thisrow{Vel_h1tS_err_M}/\LHVel},
      ]
      \addplot[fst,mono,pzero] table {\monofst};
      \addplot[fst,ac,pten]    table {\actfst};
      \addplot[fst,ac,phun]    table {\achfst};
      \addplot[scn,mono,pzero] table {\monoscn};
      \addplot[scn,ac,pten]    table {\actscn};
      \addplot[scn,ac,phun]    table {\achscn};
    \nextgroupplot[
        xlabel={Elapsed time $\times$ cores [\si{s}]},
        title={Pressure - $\ell^2(\Lebesgue2)$},
        table/x=Tot_time,table/y expr={\thisrow{Prs_l2tS_prsc}/\LLPrs},
        legend pos=south west
      ]
      \addplot[fst,mono,pzero] table {\monofst};
      \addplot[fst,ac,pten]    table {\actfst};
      \addplot[fst,ac,phun]    table {\achfst};
      \addplot[scn,mono,pzero] table {\monoscn};
      \addplot[scn,ac,pten]    table {\actscn};
      \addplot[scn,ac,phun]    table {\achscn};
    \end{groupplot}
  \end{tikzpicture}
  \caption{3D modified Taylor--Green vortex. Comparison of computational efficiency for the various first- and second-order schemes. Legend: see \cref{fig:MTGV_cvg}}
  \label{fig:MTGV_perf}
\end{figure}
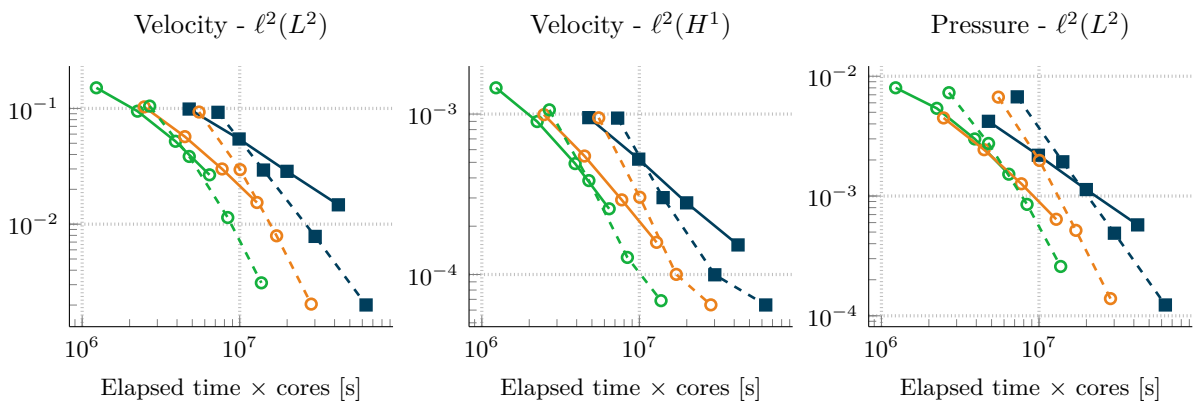

\begin{table}
  \centering
  \caption{3D modified Taylor--Green vortex. Comparison of accuracy and performance
  for the various first- and second-order schemes. Time step set to $\deltat=\frac{T}{128}$}
  \begin {tabular}{c|ccc|ccc}%
  \multicolumn {1}{c|}{}& \multicolumn {3}{c|}{$\nrm {\ErrV [\uh ]}_{\ell ^2,\cellh }$}& \multicolumn {3}{c}{$\nrm {\tGh \uh (\ErrV [\uh ])}_{\ell ^2,\cellh }$}\\Solver&First&Second&Ratio $\frac {1^{\rm {}st}}{2^{\rm {}nd}}$&First&Second&Ratio $\frac {1^{\rm {}st}}{2^{\rm {}nd}}$\\\hline %
  Monolithic&\pgfutilensuremath {2.86e{-}2}&\pgfutilensuremath {7.84e{-}3}&\pgfutilensuremath {3.6}&\pgfutilensuremath {2.80e{-}4}&\pgfutilensuremath {9.96e{-}5}&\pgfutilensuremath {2.8}\\%
  AC(10)&\pgfutilensuremath {5.20e{-}2}&\pgfutilensuremath {1.14e{-}2}&\pgfutilensuremath {4.6}&\pgfutilensuremath {4.91e{-}4}&\pgfutilensuremath {1.28e{-}4}&\pgfutilensuremath {3.8}\\%
  AC(50)&\pgfutilensuremath {3.00e{-}2}&\pgfutilensuremath {7.93e{-}3}&\pgfutilensuremath {3.8}&\pgfutilensuremath {2.92e{-}4}&\pgfutilensuremath {1.00e{-}4}&\pgfutilensuremath {2.9}\\%
  \hline \multicolumn {1}{c|}{}& \multicolumn {3}{c|}{$\nrm {\ErrP [\uh ]}_{\ell ^2,\cellh }$}& \multicolumn {3}{c}{Elapsed$\times $cores [\si {s}]}\\ Solver & First & Second & Ratio $\frac {1^{\rm {}st}}{2^{\rm {}nd}}$ & First & Second & Ratio $\frac {1^{\rm {}st}}{2^{\rm {}nd}}$\\ \hline Monolithic&\pgfutilensuremath {1.13e{-}3}&\pgfutilensuremath {4.89e{-}4}&\pgfutilensuremath {2.3}&\pgfutilensuremath {2.00e{+}7}&\pgfutilensuremath {3.01e{+}7}&\pgfutilensuremath {0.7}\\%
  AC(10)&\pgfutilensuremath {3.00e{-}3}&\pgfutilensuremath {8.51e{-}4}&\pgfutilensuremath {3.5}&\pgfutilensuremath {3.91e{+}6}&\pgfutilensuremath {8.39e{+}6}&\pgfutilensuremath {0.5}\\%
  AC(50)&\pgfutilensuremath {1.27e{-}3}&\pgfutilensuremath {5.17e{-}4}&\pgfutilensuremath {2.4}&\pgfutilensuremath {7.73e{+}6}&\pgfutilensuremath {1.71e{+}7}&\pgfutilensuremath {0.5}\\%
  \end {tabular}%
  \label{tab:MTGV}
\end{table}

We can now proceed to the assessment of the computational efficiency of the various
schemes with the tolerances prescribed as above. The performances of the three
strategies (monolithic, AC($\ACprm=10$), and AC($\ACprm=50$)), for both first- and
second-order schemes, are compared in \cref{fig:MTGV_perf}. The AC method turns
out to be more efficient: having fixed the order in time for the three compared
methods, the AC($\ACprm=10$) and AC($\ACprm=50$) methods can achieve a given error
threshold ($y$-axis) in less computation time ($x$-axis) than the monolithic method.
The performance of AC depends on the parameter $\ACprm$: the higher $\ACprm$, the
more numerical effort (hence, computational time) is required, and altogether the
choice $\ACprm=10$ turns out to be more effective than the choice $\ACprm=50$. When
considering the second-order schemes, the difference between the AC and
monolithic approaches is somewhat mitigated. One explanation is that the bootstrap
procedure employed in the second-order AC requires two linear system resolutions per
time step, whereas only one is needed with the BDF2 and monolithic coupling. Finally,
the results show an advantage of the second-order schemes over the first-order
ones, especially when accurate calculations are considered. A comparison at the
fixed time-step $\deltat=\frac{T}{128}$ (corresponding to the finest time step before
the space discretization errors tend to become influential, see the discussion of
\cref{fig:MTGV_cvg}) is presented in \cref{tab:MTGV}. Focusing first
on the bottom right part of the table related to the computational times, one
can compare the performances of monolithic and AC methods, with the 
latter being consistently two times faster than
the former, while still providing satisfactory errors. Moreover, the errors obtained with
second-order schemes are more than two times smaller than those obtained with
first-order schemes, while the computation time are less than two times larger.

\section{Conclusions}
\label{sec:conc}

We have extended the CDO-Fb schemes to the unsteady incompressible Navier--Stokes
problem and investigated the accuracy and efficiency of the artificial
compressibility (AC) technique for the time discretization using either a first-order
time-stepping scheme or a second-order one. The assessment of the AC technique is
made by means of  systematic comparisons with the standard monolithic approach, using
either a first-order or a second-order version.  First, optimal convergence rates in
time have been recovered for all the time-stepping schemes in various norms for the
velocity and the pressure. Moreover, the stability of the time-stepping schemes when
the nonlinear convection term is treated explicitly has been investigated on
Cartesian and polygonal meshes, showing in both cases a linear dependency on the
reciprocal of the Reynolds number, with only a slightly tighter restriction when
using polygonal meshes. Finally, an assessment of the performance of all the schemes
has been carried out using large 3D polytopal meshes. The AC method proved to be an
efficient alternative to the classical monolithic approach: it ensures accuracy
levels close to those obtained with the monolithic approach while using only half of
the computational time. Moreover, a comparison between first- and second-order
time-stepping schemes indicates an advantage for the latter: second-order schemes
provide errors as much as four times smaller than those obtained with first-order
schemes, while taking only up to twice the computation time. 

As an outlook, the improvement of the preconditioning for large linear systems arising from the AC
technique is paramount to achieve scalable and higher efficiency since the presence of the
$\graddiv$ term prevents the usual preconditioners such as algebraic multigrid techniques to
perform well. The construction of a specific preconditioner adapted from the seminal work in
\cite{ArFaW97} along with the more recent ones in~\cite{BeOlW11} and~\cite{FaMiW19} will be the
subject of further investigations.  

\section*{Acknowledgments}
The PhD fellowship of R. Milani was partially supported by EDF R\&D and ANRT.

\bibliographystyle{siam}
\bibliography{library}

\end{document}